\documentclass[11pt]{article}

\usepackage{amsfonts,amssymb,amsmath}
\hsize \textwidth
\advance \hsize by -\marginparwidth
\evensidemargin \oddsidemargin
\advance\hoffset by 5mm
\newcommand{\pdr}[2]{\frac{\partial{#1}}{\partial{#2}}}
\newcommand{\Rm}{{\mathbb R}}
\newcommand{\Pm}{{\mathbb P}}

\def\R{\mathbb{R}}
\newcommand{\eps}{\varepsilon}
\newcommand{\commentout}[1]{}

\def\N{\mathbb{N}}

\newcommand{\farc}{\frac}
\newcommand{\no}{\nonumber}

\newcommand{\norm}[1]{\lVert #1 \rVert}
\newcommand{\qed}{$\Box$}
\newcommand{\br}{\begin{eqnarray}}
\newcommand{\er}{\end{eqnarray}}
\newcommand{\be}{\begin{equation}}
\newcommand{\ee}{\end{equation}}
\newcommand{\baa}{\begin{array}}
\newcommand{\eaa}{\end{array}}
\newcommand{\ba}{\begin{eqnarray}}
\newcommand{\ea}{\end{eqnarray}}
\def\di{\displaystyle}
\newtheorem{theorem}{\bf Theorem}[section]
\newtheorem{theo}[theorem]{Theorem}
\newtheorem{thm}[theorem]{Theorem}

\newtheorem{lem}[theorem]{Lemma}
\newtheorem{prop}[theorem]{Proposition} 
\newtheorem{cor}[theorem]{Corollary} 
 
\begin{document}
\date{} 
\title{The logarithmic delay of KPP fronts in a periodic medium}
\author{
Fran{\c c}ois Hamel\thanks{Universit\'e d'Aix-Marseille, LATP, 39 rue F. Joliot-Curie, 13453 Marseille Cedex 13, France \& Institut Universitaire de France; francois.hamel@univ-amu.fr}
\and James Nolen\thanks{Department of Mathematics, Duke University, Durham, NC 27708, USA; nolen@math.duke.edu}
\and Jean-Michel Roquejoffre\thanks{Institut de Math\'ematiques (UMR CNRS 5219), Universit\'e Paul Sabatier, 118 route de Narbonne, 31062 Toulouse cedex, France; roque@mip.ups-tlse.fr} 
\and Lenya Ryzhik\footnote{Department of Mathematics, Stanford University, Stanford CA, 94305, USA; ryzhik@math.stanford.edu}
}

\maketitle
\begin{abstract}
We consider solutions of the KPP-type equations with a periodically varying reaction rate, and compactly supported initial data. It has been shown by Bramson~\cite{bram1,bram2} in the case of the constant reaction rate that the lag between the position of such solutions and that of the traveling waves grows as $(3/2)\log t$, as $t\to+\infty$. We generalize this result to the periodic case. 
\end{abstract}

%%%%%%%%%%%%%%%%%%%%%%%%%%%%%%%%%%%%%%%%%
%%%%%%%%%%%%%%%%%%%%%%%%%%%%%%%%%%%%%%%%%

\section{Introduction}

We study solutions $u(t,x)$ of the initial value problem
\be\label{nonlinperiodic}\left\{\baa{rcll}
&&u_t  =  u_{xx}+g(x)f(u), & t>0,\ x\in\R,\vspace{3pt}\\
&&u(0,x)  =  u_0(x).\eaa\right.
\ee
The function $f$ is of class $C^1[0,1]$, and is of KPP-type. Specifically, we assume that
\be\label{hypkpp}
f(0)=f(1)=0,\ f'(0)>0,\ f'(1)<0,\ 0<f(s)\le f'(0)s\hbox{ for all }s\in(0,1),
\ee
and that there exist $s_0\in(0,1)$, $M\ge 0$ and $\alpha>0$ such that
$$f(s)\ge f'(0)s-M\,s^{1+\alpha}\hbox{ for all }s\in[0,s_0].$$ 
We assume the function $g(x) \in C^1(\Rm)$ is $1$-periodic, and that there are two constants $g_{1,2}$ such that
\[
0<g_1\le g(x)\le g_2<+\infty.
\]
By modifying the definition of $g(x)$, we may assume without loss of generality that 
\[
f'(0) = 1.
\]  
Such equations model numerous problems
in biology and other applications, and have been extensively studied since the early papers by Fisher~\cite{Fisher} and Kolmogorov,
Petrovskii and Piskunov~\cite{kpp} -- see~\cite{Xinbook} for a recent review.

We are interested in the spreading rate for solutions of (\ref{nonlinperiodic}) with the non-negative  compactly supported 
initial conditions $u_0$ that  satisfy
$$
0\le u_0\le 1, \hbox{ and }\hbox{esssup}_{\R}u_0>0.
$$
The strong parabolic maximum principle implies that $0<u(t,x)<1$ for all $t>0$ and $x\in\R$.

%%%%%%%%%%%%%%%%%%%%%%%%%%%%%%%%%%%%%%%%%

\subsubsection*{Results in a homogeneous medium}

Let us first recall what is known when the function $g(x)$ is a constant: $g(x)\equiv 1$. Then, given any $c\ge c^*={2}$, there 
exists a traveling wave solution of 
$$u_t=u_{xx}=f(u),$$
of the form $u(t,x)=U_c(x-ct)$. The function $U_c$ satisfies
$$-cU_c'=U_c''+f(U_c),~~U_c(-\infty)=1,~~U_c(+\infty)=0,\ \ 0<U_c<1.$$
For $c > c^*$ the function $U_c(x)$ decays exponentially as $x\to+\infty$: $U_c(x)\sim Ce^{-\lambda_cx}$, with the decay rate 
$\lambda_c$ being the smallest positive solution of
$$\lambda^2-c\lambda+1=0.$$
On the other hand, at $c=c^*$ the traveling wave asymptotics is $U_{c^*}(x)\sim Cxe^{-\lambda^*x}$, with $\lambda^*=1$.
It has been shown in the pioneering work of Bramson~\cite{bram1,bram2} that solutions of the initial value problem~(\ref{nonlinperiodic})
with compactly supported initial data $u_0(x)$ ``are located" (on the right half-line $\R_+=[0,+\infty)$) at 
$$X(t)=c^*t-\frac{3}{2\lambda^*}\log t+O(1)\hbox{ as $t\to +\infty$. }$$
More precisely, $u(t,x)$ satisfies the following property: given any $\eps>0$ if we set
\[
X_\eps(t)=\sup\big\{x\in\Rm:~u(t,x)\ge\eps\big\},~~Y_\eps(t)=\inf\big\{x\in\Rm_+:~u(t,x)\le 1-\eps\big\},
\]
then
$$X_\eps(t)=c^*t-\frac{3}{2\lambda^*}\log t+O(1)\hbox{ as $t\to +\infty$, }$$
and
$$Y_\eps(t)=c^*t-\frac{3}{2\lambda^*}\log t+O(1)\hbox{ as $t\to +\infty$. }$$
In other words, the region in $\R_+$ where $u(t,x)$ transitions from the value $u\approx 1$ to $u\approx 0$ has a width that is uniformly
bounded in time, and is located at the distance $(3/2\lambda^*)\log t$ behind the location of the traveling wave with minimal speed $c^*$. Bramson's proofs
were based on probabilistic techniques, and were later extended by G\"artner to higher dimensions~\cite{Gartner}, and recently
revisited by Roberts~\cite{Roberts},
 while
a PDE proof of this result was later given by \cite{Lau} with the additional assumption $f'(s)\le f'(0)$ on $[0,1]$, and 
recently in the companion paper \cite{HNRR1}, with other results in this direction
obtained earlier in~\cite{Uchiyama}.

We should also mention a very interesting paper~\cite{FZ} where the medium is taken to be time-dependent, with the reaction coefficient
taking two different values $\sigma_1$ and $\sigma_2$ on the time intervals $[0,T]$ and $[T,2T]$. It is shown by probabilistic techniques 
that the lag behind $X(t)$ and traveling front position  depends strongly on whether $\sigma_1>\sigma_2$ or
$\sigma_2>\sigma_1$.

%%%%%%%%%%%%%%%%%%%%%%%%%%%%%%%%%%%%%%%%%

\subsubsection*{Periodic pulsating fronts}

In order to understand how Bramson's results can be adapted to a periodic environment, let us recall the notion of a pulsating
traveling wave that generalizes the notion of a traveling wave to periodic media. 
A pulsating front with speed $c > 0$ is a function $U_c(t,x)$ satisfying
\be
U_t = U_{xx} + g(x)f(U), \quad x \in \Rm, \;\; t \in \Rm, \label{pulsedUeqn}
\ee
and 
\[
U(t + \frac{1}{c},x) = U(t,x - 1),
\]
as well as the boundary conditions $U(t,-\infty) = 1$, $U(t,+\infty) = 0$. Let us now recall some of the results about spreading speeds and pulsating traveling waves $U_c(t,x)$~\cite{BH,BHN,Hamel1,hr,w,Xinbook}. It is known that there is a minimal speed $c^* > 0$ such that for each $c \geq c^*$, there exists a unique up to time-shifts pulsating traveling front $U_c(t,x)$, while no pulsating traveling front exists with a speed less than $c^*$. Furthermore, all pulsating traveling fronts are necessarily increasing in $t$. Lastly, the minimal speed $c^*$ may be characterized as follows. Given $\lambda>0$, let $\psi = \psi(x,\lambda) > 0$ be the principal eigenfunction of the $1$-periodic eigenvalue problem
\begin{equation}
\psi_{xx} - 2 \lambda \psi_x + (\lambda^2 + g(x)f'(0)) \psi=\gamma(\lambda) \psi ,
~~\psi(x+1,\lambda)=\psi(x,\lambda),~~\psi(x,\lambda)>0,\quad x \in \Rm, \label{evpperiodic}
\end{equation}
and $\gamma(\lambda)$ the corresponding eigenvalue. 
The eigenfunction is normalized so that
\begin{equation}\label{feb2312}
\int_0^1 \psi(x,\lambda) \,dx = 1,
\end{equation}
for all $\lambda > 0$.  
The minimal wave speed is given by
\[
c^* = \min_{\lambda > 0} \farc{\gamma(\lambda)}{\lambda} = c(\lambda^*).
\] 
Here $\lambda^*>0$ minimizes $\gamma(\lambda)/\lambda$. In particular, we have 
\begin{equation}\label{0523-14}
\gamma'(\lambda^*) =\frac{\gamma(\lambda^*)}{\lambda^*} = c^*.
\end{equation}

%%%%%%%%%%%%%%%%%%%%%%%%%%%%%%%%%%%%%%%%%

\subsubsection*{The main results}

Our first main result is as follows. 

\begin{thm}\label{thm-delay-apr5}
Let $u(t,x)$ be a solution of (\ref{nonlinperiodic}) with the initial data $u_0(x)$ such that $0\le u_0(x)\le 1$, $u_0(x)\not\equiv 0$,
and $u_0(x)=0$ for $|x|>M$ with some $M>0$. Then for any $\eps>0$ there exist 
$s(\eps)$ and $L(\eps)$
so that 
$$u(t,x)\ge 1-\eps\hbox{ for all }t>s(\eps)\hbox{ and all }x\in\big[0 \, , \,c^*t-\frac{3}{2\lambda^*}\log t - L(\eps)\big]$$
and
$$u(t,x) < \eps\hbox{ for all }t > s(\eps)\hbox{ and all }x\in\big[c^*t-\frac{3}{2\lambda^*}\log t + L(\eps) \, , \, +\infty\big).$$
\end{thm}

This generalizes directly Bramson's results to a periodic medium: the front is located at distance $(3/2\lambda^*)\log t$ behind the
pulsating front. 

Let us explain informally how the logarithmic decay comes about.  The main observation, from the PDE point of view, is that
solutions of the nonlinear problem (\ref{nonlinperiodic}) behave very similar to those of the linearized problem
$$v_t=v_{xx}+g(x)v,$$
with the Dirichlet boundary condition $v(t,c^*t)=0$. In the homogeneous case, with $g(x)\equiv 1$, $c^*=2$ and $\lambda^*=1$, 
let us write 
$v(t,x)=p(t,x)e^{-(x-2t)}$.
Then $p(t,x)$ satisfies
\[
p_t=p_{xx}-2p_x,~~x>2t,~~p(t,2t)=0.
\]
Changing variables $y=x-2t$, we get
\[
p_t=p_{yy},~~p(t,0)=0.
\]
It follows that $p(t,y=1)\sim t^{-3/2}$ as $t\to +\infty$, or, in the original variables $v(t,x=2t+1)\sim t^{-3/2}$.  
Assuming that the solution $u(t,x)$ of the nonlinear problem has the same  behavior as $v(t,x)$, and has the 
exponential asymptotics
$u(t,x)\sim e^{-(x-X(t))}$, we deduce that $X(t)\sim 2t-(3/2)\log t$. For the homogeneous case $g \equiv 1$, we have worked out this argument in detail in \cite{HNRR1}. The bulk of the proof in the periodic case is in getting the decay estimates for the heat kernel
in a half space with periodic coefficients. These estimates are well known in the whole space~\cite{FS,Nash,Norris} but we are not aware
of such results in a half space for periodic coefficients.

In the proof of Theorem~\ref{thm-delay-apr5}, one shows actually more precise exponential estimates on $u(t,x)$ for $x\ge c^*t-(3/(2\lambda^*))\log t$. These estimates imply that the solution~$u$ is asymptotically trapped between two finite space-shifts of the minimal front $U_{c^*}$ around the position $x=c^*t-(3/(2\lambda^*))\log t$. Equivalently, $u$ is asymptotically trapped between two finite time-shifts of the minimal front $U_{c^*}$ around the time $t-(3/(2c^*\lambda^*))\log t$. Then, by passing to the limit along any level set, any limiting solution is necessarily equal to a shift of the minimal front: this follows from a new Liouville-type result which is similar to what had already been known in the homogeneous case. For more details, we refer to Section 8, where the following result is proved:

\begin{theo}\label{th2}
There exist a constant $C\ge 0$ and a function $\xi:(0,+\infty)\to\R$ such that $|\xi(t)|\le C$ for all $t>0$ and
\be\label{fronts}
\lim_{t \to +\infty} \,\left \lVert u(t,\cdot) - U_{c^*}\big(t-\frac{3}{2c^*\lambda^*}\log t+\xi(t),\cdot\big) \right \rVert_{L^\infty(0,+\infty)}  = 0.
\ee
Furthermore, for every $m\in(0,1)$ and every sequence $(t_n,x_n)$ such that $t_n\to+\infty$ and $x_n-[x_n]\to x_{\infty}\in[0,1]$ as $n\to+\infty$, and $u(t_n,x_n)=m$ for all $n\in\N$, there holds
\be\label{fronts2}
u(t+t_n,x+[x_n])\mathop{\longrightarrow}_{n\to+\infty} U_{c^*}(t+T,x)\hbox{ locally uniformly in }(t,x)\in\R^2,
\ee
where $[x_n]$ denotes the integer part of $x_n$ and $T\in\R$ denotes the unique real number such that $U_{c^*}(T,x_{\infty})=m$.
\end{theo}

Theorem~\ref{th2} shows in particular the convergence to the family of shifted  minimal fronts along the level sets of the solution $u$. Results of this type have been obtained recently in~\cite{dgm} for more general nonlinearities $f$ and Heaviside initial conditions $u_0$ and in~\cite{Giletti} for asymptotically periodic KPP functions~$f$ and compactly supported initial conditions~$u_0$. The proofs in~\cite{dgm,Giletti} are completely different from the ones used here: they are based on the time-decay property of the number of intersections of any two solutions and on the fact that the minimal fronts are the steepest ones. They hold for more general functions $f$ but do not provide the logarithmic shift of the position of the solutions.

%%%%%%%%%%%%%%%%%%%%%%%%%%%%%%%%%%%%%%%%%

\subsection*{Connection to branching Brownian motion}

When $g$ is constant and $f(u) = u(1 - u)$, there is a well-known connection between solutions 
of~(\ref{nonlinperiodic}) and branching Brownian motion \cite{bram1,mck}. 
Consider a branching Brownian motion with constant branching rate $g > 0$. Initially, there is one Brownian particle, $X_1(0) = 0$. 
At a random time $T_1$, which is an independent exponential random variable with rate $g$, this particle gives birth to two 
independent Brownian motions and then dies immediately itself. 
The two new particles start their motions from the final location of the parent particle. The process continues in this way, each living particle reproducing and dying at an independent random time, leaving two new Brownian particles as offspring.  As shown by McKean \cite{mck}, the function 
\be
u(t,x) = \Pm \left( \max_{k \in L(t) } X_k(t) > x \;|\; X_1(0) = 0\right) \label{constbranch}
\ee
satisfies
\[
u_t = \frac{1}{2} u_{xx} + g u(1 - u).
\]
and 
\[
u(0,x) = 1, \quad x \leq 0; \quad \quad u(0,x) = 0, \quad x > 0.
\]
The set $L(t)$ in (\ref{constbranch}) denotes the set of indices corresponding to particles that are alive at time $t$. 

When $g(x)$ is not constant, there is a similar interpretation of (\ref{nonlinperiodic}) in terms of a branching Brownian motion with space-dependent branching rate $g(x) > 0$. In that case, we start the particle initially at $x$: $X_1(0) = x$. At a random time, this particle produces two Brownian offsprings and then dies immediately. The branching time is constructed from an independent exponential random variable: if $S$ is a standard exponential random variable, independent of $X_1(t)$, then the time at which $X_1$ branches is
\[
T_1 = \inf \left\{ t > 0 \;|\; S < \int_0^t g(X_1(r))\,dr \right\},
\]
so that $T_1$ satisfies $\Pm(T_1 > t \;|\; X_1) = e^{- \int_0^t g(X_1(r)) \,dr}$.  Using arguments as in \cite{bram1,mck}, 
one can show that the function
\be
u(t,x) = \Pm \left( \min_{k \in L(t)} X_k(t) < 0 \;|\; X_1(0) = x\right) \label{gxbranch}
\ee
satisfies
\[
u_t = \frac{1}{2} u_{xx} + g(x) u(1 - u).
\]
and 
\[
u(0,x) = 1, \quad x \leq 0; \quad \quad u(0,x) = 0, \quad x > 0.
\]
When $g$ is constant, it is easy to see that the two formulas (\ref{constbranch}) and (\ref{gxbranch}) define the same function. However, if $g$ is not constant then (\ref{gxbranch}) need not be equivalent to (\ref{constbranch}).

The zero Dirichlet boundary condition exactly corresponds to G\"artner's~\cite{Gartner} strategy of killing the branching Brownian motion 
at a moving boundary. 

The paper is organized as follows. Section~\ref{sec:lower} contains the basic elements of the 
proof of the lower bound for the solution, while the main steps of the proof of the
upper bound are contained in Section~\ref{sec:upper}. Sections 4, 5, 6 and 7 contain the proofs of the auxiliary
results formulated in these two sections. Lastly, Section 8 is devoted to the proof of Theorem~\ref{th2}.\hfill\break

{\bf Acknowledgment.}  This work was motivated by a series of lectures given by Eric Brunet at Banff Conference Center in March 2010
on his work with Bernard Derrida concerning KPP and related particle models.
JN was supported by NSF grant DMS-1007572, and LR by NSF grant DMS-0908507. FH and JMR were supported by ANR grant PREFERED.

%%%%%%%%%%%%%%%%%%%%%%%%%%%%%%%%%%%%%%%%%
%%%%%%%%%%%%%%%%%%%%%%%%%%%%%%%%%%%%%%%%%

\section{The lower bound: outline of the proof}\label{sec:lower}

%%%%%%%%%%%%%%%%%%%%%%%%%%%%%%%%%%%%%%%%%

\subsubsection*{The linearized Dirichlet problem}

The proof of the lower bound in Theorem \ref{thm-delay-apr5} is based on the analysis of the 
linearized problem with the Dirichlet boundary condition at $x = c^*t$ (recall that $f'(0)=1$):
\begin{eqnarray}\label{2803-6}
&&w_t=w_{xx}+g(x)w,~~~~~x \geq c^*t,\\
&& w(t,c^*t)=0, \quad t \geq 0,\nonumber\\
&&w(0,x)=u_0(x), \quad x \geq 0.\nonumber
\end{eqnarray}
As we will see, with an appropriate choice of $a(t)$, the function $\bar w(t,x) = a(t) w(t,x)$ will be a subsolution of the 
nonlinear equation (\ref{nonlinperiodic}). Therefore, a lower bound on $u$ will follow from a lower bound on $w$.  
It is convenient to represent $w(t,x)$ in the form
\begin{equation}\label{feb2322}
w(t,x)=e^{-\lambda^*(x-c^*t)}\psi(x,\lambda^*)p(t,x).
\end{equation}
Here $\psi(x,\lambda^*)$ is the eigenfunction of (\ref{evpperiodic})-(\ref{feb2312}),
with $\lambda=\lambda^*$ satisfying~(\ref{0523-14}), and  $p(t,x)$ satisfies 
\begin{eqnarray}\label{2803-8bis}
&&p_t=p_{xx}+\farc{2\phi_x}{\phi} p_x, \quad x \geq c^*t,\\
&& p(t,c^*t) = 0, \quad t > 0\nonumber \\
&&p(0,x)=p_0(x) = u_0(x) e^{\lambda^* x} (\psi(x,\lambda^*))^{-1}, \quad x > 0,\nonumber
\end{eqnarray}
with $\phi(t,x) = e^{-\lambda^*(x - c^* t)} \psi(x,\lambda^*)$.
The initial data $p_0(x)$ is nonnegative and compactly supported on $[0,+\infty)$.  For convenience, we define the function
\begin{equation}\label{0523-28}
\kappa(x)=\frac{2\phi_x}{\phi}=-2\lambda^*+2\frac{\psi_x(x,\lambda^*)}{\psi(x,\lambda^*)},
\end{equation}
which is the drift term   in (\ref{2803-8bis}). This function $\kappa(x)$ is $1$-periodic in $x$, and is independent of $t$. 

The first (and longest) step in the proof of the lower bound in
Theorem \ref{thm-delay-apr5} is the following lower bound on $p(t,x)$, which implies a lower bound on $w(t,x)$.  

\begin{prop}\label{cor-0524}
There exist constants $T_0 > 0$, $\sigma > 0$, and $C_0 > 0$ such that
$$\displaystyle p(t,c^*t+\sigma \sqrt{t})\ge\farc{C_0}{t}\ \hbox{ for all }t \geq T_0.$$
\end{prop}

For the homogeneous medium, when $g$ is constant, it is rather simple to derive the bound in Proposition \ref{cor-0524}. 
In that case $\psi(x)\equiv 1$, and $\phi = e^{-\lambda^* (x - c^* t)}$, so that $\kappa \equiv - 2 \lambda^*$. 
Moreover, when $g$ is constant it happens that $2 \lambda^* = c^*$, so the function $z(t,x)= p(t,x + c^*t)$ 
satisfies the heat equation $z_t = z_{xx}$ on the half-line with Dirichlet boundary condition $z(t,0) = 0$. Then, using the explicit formula one finds that there exists $C > 0$ so that
\be
\farc{x-c^*t}{Ct^{3/2}}  \leq p(t,x) \leq \farc{C(x-c^*t)}{t^{3/2}} \label{pexplheat}
\ee
holds for $x \in [c^*t, c^*t + \sqrt{t}]$. When $g$ is not constant, however, the analysis is more difficult: it is not generally true that $2 \lambda^* = c^*$, nor do we have an explicit formula for the heat kernel associated with~(\ref{2803-8bis}). Moreover, the standard bounds for the heat kernel for equation (\ref{2803-8bis}) on the entire line $x \in \Rm$ do not immediately imply the needed estimate for the Dirichlet problem on the half-line $x > c^* t$.

%%%%%%%%%%%%%%%%%%%%%%%%%%%%%%%%%%%%%%%%%

\subsubsection*{From the linearized problem to a subsolution for the nonlinear problem}

Given the lower bound of Proposition \ref{cor-0524}, the next step is to construct a subsolution for (\ref{nonlinperiodic}) using the solution of (\ref{2803-6}). If $\bar w(t,x)=a(t)w(t,x)$, then $\bar w(t,x)$ is a subsolution for (\ref{nonlinperiodic}),
that is,
\[
\bar w_t\le\bar w_{xx}+g(x)\bar w-g(x)q(\bar w),
\]
with $q(\bar w)=\bar w-f(\bar w)$, provided that
\begin{equation}\label{0524-02}
a'(t)\,w(t,x)\le -g(x)q(a(t)w(t,x)).
\end{equation}
As $q(s)\le ms^2$, and $g(x)$ is uniformly bounded from above and below by two positive constants,~(\ref{0524-02}) holds provided that
\begin{equation}\label{0524-04}
a'(t)\,w(t,x)\le -Ma(t)^2w(t,x)^2,
\end{equation}
with a large enough constant $M$. We claim that there exists a constant $C_0>0$, depending on the initial data $u_0$, such that 
\begin{equation}\label{0524-06}
w(t,x)\le \frac{C_0}{(t+1)^{3/2}}.
\end{equation}
for all $t \geq 0$ and $x \in \Rm$ (we may define $w(t,x) = 0$ for $x < c^*t$).  This estimate is a consequence of an upper bound on $p(t,x)$: 

\begin{lem}\label{prop-0504-2bis}
There exists a constant $C>0$ such that
\begin{equation}\label{3103-02bis}
|p(t,x+c^*t)|\le \frac{Cx}{(t+1)^{3/2}}\int_{0}^\infty yp_0(y)\,dy
\end{equation}
for all $t > 0$ and $x > 0$.
\end{lem}

Once again, in the homogeneous case, (\ref{3103-02bis}) follows trivially from the explicit solution formula. With (\ref{3103-02bis}) in hand, using (\ref{feb2322}), we have
\[
\sup_{x \geq c^*t} w(t,x) \leq \norm{\psi(\cdot,\lambda^*)}_\infty \frac{C}{(t+1)^{3/2}} 
\left( \int_{0}^\infty yp_0(y)\,dy \right) \sup_{x \geq c^*t}\left[ e^{-\lambda^*(x - c^*t)} (x- c^*t)\right],
\]
which implies   (\ref{0524-06}). 
Next, given  (\ref{0524-06}), (\ref{0524-04}) holds provided that
$$a'(t) \le -\farc{M }{(t+1)^{3/2}}a(t)^2,$$
and we may take
\[
a(t)=\frac{a(0)}{1 + 2M a(0)(1 - (t + 1)^{-1/2})},\ \ a(0)>0,
\]
which satisfies 
\[
\frac{a(0)}{1 + 2Ma(0)} \leq a(t) \leq a(0)
\]
for all $t \geq 0$. If $a(0) < 1$, then $\bar w(0,x) \leq u_0(x)$ for all $x \in \Rm$. Therefore, the comparison principle implies
\[
u(t,x) \geq \bar w(t,x) = a(t) w(t,x) \ge C w(t,x)\ \hbox{ for all }t\ge 0\hbox{ and }x\ge c^*t.
\]
In particular, Proposition \ref{cor-0524} implies that
\begin{equation} \label{0524-10}
u(t,ct + \sigma \sqrt{t} ) \geq C t^{-1} e^{-\lambda^* \sigma\sqrt{t} } 
\end{equation}
for $t \geq T_0$.

%%%%%%%%%%%%%%%%%%%%%%%%%%%%%%%%%%%%%%%%%

\subsubsection*{From a lower bound on the far right to the bound at the front}

Now we show that (\ref{0524-10}) (a bound far on the right) implies the lower bound in Theorem~\ref{thm-delay-apr5}.  
Let $\eps > 0$. We will use (\ref{0524-10}) to show that there is a constant $L(\eps) \in \Rm$ such that
\begin{equation}\label{080402}
u(t,x) \geq 1 - 2 \eps , \quad \forall \;\; x \in \big[0,c^*t-\frac{3}{2\lambda^*}\log t - L(\eps)\big].
\end{equation}
if $t$ is sufficiently large.

Let~$k$ be a~$C^1[0,1-\varepsilon]$ function such that $k \le f$ in $[0,1-\varepsilon]$, $k(0)=k(1-\varepsilon)=0$, $k'(0)=f'(0)=1$, 
and $k>0$ on $(0,1 - \varepsilon)$. The function~$k$ then satisfies
$$k(s)\le f(s)\le f'(0)s=k'(0)s\hbox{ for all }s\in[0,1 - \varepsilon].$$
Thus, there exists a pulsating traveling front $U^k_{c^*}(t,x)$ solution of~(\ref{pulsedUeqn}) 
with nonlinearity~$k$ instead of~$f$, having the same minimal speed $c^*$, and such that $0 < U^k_{c^*} < 1 - \varepsilon$, and 
\be
\lim_{x \to -\infty} U^k_{c^*}(t,c^*t + x) = 1 - \varepsilon, \quad \lim_{x \to +\infty} U^k_{c^*}(t,c^*t + x)=0, \label{Uklim}
\ee
uniformly in $t$. Moreover, $U^k_{c^*}$ is monotone increasing in $t$.

To show (\ref{080402}), we will bound $u$ from below by the function
\be\label{deftildeU}
\tilde U(t,x)=U^k_{c^*}(t-r(t),x).
\ee
Since we have $\partial_t U^k_{c^*}(t,x) > 0$ for all $t$ and $x$, the function $\tilde U(t,x)$ satisfies
$$\tilde U_t-\tilde U_{xx}-g(x) k(\tilde U)= (1 - r'(t)) \partial_t U^k_{c^*} - \partial_x^2 U^k_{c^*} -g(x) k(U^k_{c^*})=-r'(t) \partial_t U^k_{c^*}\le 0,$$
provided that $r'(t)\ge 0$.  In this case, since $f\ge k$ in $[0,1-\eps]$, $\tilde U$ is a subsolution of the equation 
\[
u_t - u_{xx} - g(x)f(u) = 0.
\]
Since $g(x) \geq g_1 > 0$, it is known from~\cite{aw} that $u(t,x)\to 1$ as $t\to+\infty$ locally uniformly in $x\in\R$. 
Therefore, there exists $T_1 > 0$, depending on $u_0$ and $\varepsilon$, such that $u(t,0)\ge 1 - \varepsilon$ for all $t\ge T_1$. Therefore, 
\[
\tilde U(t,0) < 1 - \varepsilon \le u(t,0), \quad \forall t \geq T_1.
\]
By taking $T_1$ larger, if necessary, we may assume $T_1 > T_0$ so that (\ref{0524-10}) holds for all $t \geq T_1$. Therefore, the maximum principle and (\ref{0524-10}) imply that the bound
\begin{equation}\label{080414}
\tilde U(t,x)\le u(t,x) \quad \text{for all} \quad x \in [0, c^*t+ \sigma \sqrt{t}], \;\; t \geq T_1,
\end{equation}
will hold, if both
\begin{equation}
\tilde U(T_1,x)  \le u(T_1,x), \quad x \in [0, c^*T_1 + \sigma \sqrt{T_1}] \label{UT0}
\end{equation}
and 
\begin{equation}\label{080408}
\tilde U(t,c^*t + \sigma \sqrt{t})\le \farc{C}{t} e^{-\lambda^* \sigma \sqrt{t}}, \quad t > T_1 
\end{equation}
are satisfied. 

Let us now verify that (\ref{UT0}) and (\ref{080408}) hold with
\be\label{defr}
r(t)=\left(\farc{3}{2\lambda^* c^*}\right)\log t +  L_0,
\ee
if $L_0$ is sufficiently large. Because (\ref{Uklim}) holds uniformly in $t$, it is clear that for $T_1$ fixed, we may take $L_0$ sufficiently large (depending only on $u_0$ and $T_1$) so that (\ref{UT0}) is satisfied. It was shown in~\cite{Hamel1} that 
the function $U^k_{c^*}(t,x)$ satisfies
$$U^k_{c^*}(t,x)\le C(x-c^*t)e^{-\lambda^*(x-c^*t)},$$
for $x\ge c^*t + 1$. Hence,  the function $\tilde U$ satisfies
$$\begin{array}{rcl}
\tilde U(t,c^*t+\sigma \sqrt{t}) & = & U^k_{c^*}(t-r(t),c^*t+\sigma \sqrt{t})\\
& \le & C(c^*t+\sigma \sqrt{t}-c^*(t-r(t)))e^{-\lambda^*(c^*t+\sigma \sqrt{t}-c^*(t-r(t)))}\nonumber\\
& = & C \left(\sigma \sqrt{t} + \displaystyle\frac{3}{2\lambda^*} \log t  + c^* L_0\right)\displaystyle\frac{e^{-\lambda^* \sigma \sqrt{t}}}{t^{3/2} }e^{-\lambda^*c^*L_0}
\leq  \displaystyle\frac{C}{t} e^{-\lambda^* \sigma \sqrt{t}},
\nonumber\end{array}$$
for all $t > T_1$, provided that $L_0$ and $T_1$ are sufficiently large. Hence, (\ref{080408}) also holds for large enough $L_0$.  Therefore, (\ref{080414}) must hold for large enough $L_0$ and $T_1$. 

For $t \geq T_1$ and $h > 0$, let $A_t$ be the interval
\[
A_t =\big[0, c^*t - \frac{3}{2\lambda^*} \log(t) -  c^* L_0 - h\big] \subset [0, c^* t + \sigma \sqrt{t}].
\] 
We have now shown that for $t \geq T_1$,
\begin{eqnarray} \label{infinfU}
 \inf_{x \in A_t} u(t,x) & \geq & \inf_{x \in A_t} U^k_{c^*} \left(t - \frac{3}{2c^*\lambda^*}\log(t) - L_0, x \right) \geq  \inf_{s \in \Rm} \inf_{y \leq c^* s - h} U^k_{c^*} \left(s, y \right) 
\end{eqnarray}
From the properties (\ref{Uklim}) of $U^k_{c^*}$, we know that the right side of (\ref{infinfU}) is larger than $(1 - 2\varepsilon)$ if $h > 0$ is sufficiently large. This proves (\ref{080402}). 

Thus, we have reduced the proof of the lower bound in
Theorem~\ref{thm-delay-apr5} to the proof of Proposition~\ref{cor-0524} and Lemma~\ref{prop-0504-2bis}. 
We postpone them until later sections, and first describe in Section~\ref{sec:upper} how the upper bound in this theorem is proved.

%%%%%%%%%%%%%%%%%%%%%%%%%%%%%%%%%%%%%%%%%
%%%%%%%%%%%%%%%%%%%%%%%%%%%%%%%%%%%%%%%%%

\section{The upper bound: outline of the proof}\label{sec:upper}

%%%%%%%%%%%%%%%%%%%%%%%%%%%%%%%%%%%%%%%%%

\subsubsection*{The linearized problem in the logarithmically shifted reference frame}

As we have seen, the idea behind the $(3/{2\lambda^*}) \log(t)$ delay is that the evolution is driven by the behavior 
of solutions to the Dirichlet problem (\ref{2803-6}), which is 
$$z_t-z_{xx} - g(x)z=0,~~x>c^* t,$$
with $z(t,c^*t)=0$. The problem is that such solutions that are initially compactly supported
will decay in time like $t^{-3/2}$, hence they can not serve as super-solutions to the non-linear problem.
The correction to this inconvenience is to devise a reference frame  in which the
Dirichlet problem will have solutions that remain bounded both from above and below by positive constants for finite $x$, and this is exactly what the 
$3/({2\lambda^*})\log t$ shift achieves.  We expect the front to be at $x(t)=c^*t-r\log t$, with $r=3/(2\lambda^*)$. 
For the moment, let us assume that the constant $r$ is still general, and we will choose $r$ appropriately later. Accordingly, we consider the Dirichlet problem
$$\left\{\begin{array}{l}
z_t-z_{xx}-g(x)z=0,~~t>0,~~x>c^*t-r\log (t + T) + r \log(T),\vspace{3pt}\\
z(t,c^*t-r\log (t+T) + r\log(T))=0,\nonumber\end{array}\right.$$
with a given nonnegative continuous compactly supported initial condition $z(0,\cdot)\not\equiv 0$ in $(0,+\infty)$.

Define the new time variable $\tau$ by
$$c^*\tau =c^*t-r\log (t+T)+r\log T,$$
and set $\tilde z(\tau,x)=z(t,x)$.  
Let us also denote $t=h(\tau)$, and choose $T> 0$ sufficiently large so that the function $h(\tau)$ is well defined and monotonic.
Then we have
\[
\tilde z_\tau=h'(\tau)z_t=h'(\tau)[z_{xx}+g(x)z]=h'(\tau)[\tilde z_{xx}+g(x) \tilde z],~~x>c^*\tau,
\]
and $\tilde z(\tau,c^*\tau)=0$. Next, set
\[
\tilde z(\tau,x)=e^{-\lambda^*(x-c^*\tau)} \psi(x,\lambda^*) \alpha(\tau) \tilde p(\tau,x),
\]
with an increasing function $\alpha(\tau) > 0$ to be determined. Here, as before, $\psi(x,\lambda^*)$ is the eigenfunction of 
(\ref{evpperiodic})-(\ref{feb2312}).
The function $\tilde p(\tau,x)$ must satisfy
\begin{equation}\label{oct1206}
\farc{1}{h'(\tau)} \tilde p_\tau = \tilde p_{xx} + 2\frac{\phi_x}{\phi} \tilde p_x + \left(-\frac{1}{h'(\tau)} \frac{\alpha'(\tau)}{\alpha(\tau)}+  \lambda^* c^* \left( 1 - \frac{1}{h'(\tau)} \right) \right) \tilde p = 0,~~~\tau>0,~x>c^*\tau,
\end{equation}
where $2\phi_x/\phi$ is as in~(\ref{0523-28}). We first compute $h'(\tau)$:
\[
\farc{1}{h'(\tau)}=1-\farc{r}{c^*(h(\tau)+T)}=1-\farc{r}{c^*(\tau+T)+{r}\log((t+T)/T)}=1-\farc{r}{c^*(\tau+T)}+\beta(\tau),
\]
with
\[
\beta(\tau)=\farc{r}{c^*(\tau+T)}-\farc{r}{c^*(\tau+T)+{r}\log((t+T)/T)}=
\frac{r^2\log((t+T)/T)}{c^*(\tau+T)(c^*(\tau+T)+{r}\log((t+T)/T))}.
\]
Observe that $|\beta(\tau)| \leq C |\tau^{-3/2}|$, and if $r> 0$, then $h'(\tau) > 1$ for all $\tau > 0$.

To eliminate the low-order term in (\ref{oct1206}), we now choose $\alpha(\tau)$ so that
\[
\frac{\alpha'(\tau)}{\alpha(\tau)} = c^* \lambda^* (h'(\tau) - 1) = \frac{r\lambda^*}{(\tau + T)} + O\left(\frac{1}{(\tau + T)^{3/2}}\right),
\]
hence
\begin{equation}\label{oct1706}
\alpha(\tau)=\exp[ {r\lambda^*}\ln(\tau+T) +O(\tau^{-1/2})]=(\tau+T)^{r\lambda^*}(1+O(\tau^{-1/2})).
\end{equation}
The function $\tilde p(\tau,x)$ then satisfies
\begin{eqnarray}\label{oct1220}
&&\frac{1}{h'(\tau)} \tilde p_\tau = \tilde p_{xx} + 2 \frac{\phi_x}{\phi} \tilde p_x, ~~~\tau >0,~~x>c^*\tau,
\end{eqnarray}
with the Dirichlet condition $\tilde p(\tau,c^* \tau) = 0$. Observe that if $r = 0$ (taking no logarithmic shift), and 
$h' \equiv 1$, this is identical to equation (\ref{2803-8bis}) which is satisfied by $p(t,x)$ that was used
in the construction of a sub-solution. However, we can not take $r=0$ and use $p(t,x)$ for a super-solution since 
$p(t,x)$ decays as $t^{-3/2}$ as $t\to+\infty$ while for a super-solution 
we need $p(t,x)$ to stay bounded from above and below for finite values of $x$.

To bound the function
\[
z(t,x) = \tilde z (\tau,x) = e^{-\lambda^* (x - c^* \tau)} \psi(x,\lambda^*) \alpha(\tau) \tilde p(\tau,x),
\]
we need an estimate on $\tilde p(\tau,x)$ from above and below. The main technical step in the proof of  the upper bound in
Theorem~\ref{thm-delay-apr5} is the following estimate on $\tilde p(\tau,x)$, which implies that $\tilde p$ has the same 
leading order behavior as $p$, even though $h'(\tau) \neq 1$ in   (\ref{oct1220}).
Let us set
$$\omega(\tau)=1 -  \frac{1}{h'(\tau)} = \farc{r}{c^*(\tau+T)}-\beta(\tau).$$
Observe that $\omega(\tau) \sim {r}/{c^* \tau}$ as $\tau \to \infty$, and  $|\omega(\tau)|\le{C}/{\tau}$,  $|\omega'(\tau)|\le C/\tau^2$ for $\tau>\tau_0$.  

\begin{prop}\label{thm-plan1}
Let $\tilde p(\tau,x)$ satisfy 
\begin{equation}\label{plan-eq2}
\left(1-\omega(\tau)\right)\tilde p_\tau = \tilde p_{xx} + 2 \frac{\phi_x}{\phi} \tilde p_x,~~x\ge c^*\tau,
\end{equation}
with the Dirichlet boundary condition $\tilde p(\tau,c^*\tau)=0$.  Then 
there exist constants $k,K,\tau_0>0$ so that
$$\farc{k(x-c^*\tau)}{\tau^{3/2}}\le \tilde p(\tau,x)\le\farc{K(x-c^*\tau)}{\tau^{3/2}},$$
for all $x\in(c^*\tau,c^*\tau +k\sqrt{\tau})$ and all $\tau>\tau_0$.
\end{prop}

%%%%%%%%%%%%%%%%%%%%%%%%%%%%%%%%%%%%%%%%%

\subsubsection*{Proof of the upper bound in Theorem~\ref{thm-delay-apr5}}

In terms of the function $\tilde z(\tau,x)$, Proposition \ref{thm-plan1} says that 
$$\frac{\alpha(\tau)}{\tau^{3/2}} k(x-c^*\tau)e^{-\lambda^*(x-c^*\tau)}  \le \tilde z(\tau,x)\le \frac{\alpha(\tau)}{\tau^{3/2}} {K(x-c^*\tau)}e^{-\lambda^*(x-c^*\tau)}$$
holds for all $x\in(c^*\tau,c^*\tau +k\sqrt{\tau})$ and all $\tau>\tau_0$, even if it means changing the positive constants $k$ and $K$. Expression (\ref{oct1706}) for $\alpha(\tau)$ shows that the choice of $r={3}/{(2\lambda^*)}$ gives
$$K_1\le \frac{\alpha(\tau)}{\tau^{3/2}}\le K_2, \quad \tau \geq \tau_0,$$
and therefore
$$k(x-c^*\tau)e^{-\lambda^*(x-c^*\tau)}  \le \tilde z(\tau,x)\le{K(x-c^*\tau)}e^{-\lambda^*(x-c^*\tau)}$$
holds for all $x\in(c^*\tau,c^*\tau +k\sqrt{\tau})$ and all $\tau>\tau_0$. 

Now, we go back to the $t$ variable and bound $z(t,x) = \tilde z(\tau,x)$. Since
$$c^*\tau =c^*t-r\log (t+T)+r\log T,$$
we get the lower and upper bounds
\be\label{plan-eq14}\begin{array}{l}
z(t,x)  \geq  k(x-c^*t+r\log (t+T)-r\log T)e^{-\lambda^*(x-c^*t+r\log (t+T)-r\log T)},  \\
z(t,x)  \leq  {K(x-c^*t+r\log (t+T)-r\log T)}e^{-\lambda^*(x-c^*t+r\log (t+T)-r\log T)},\end{array}
\ee
for all $t \geq h(\tau_0)$, in the interval
$$c^*t-r\log (t+T)+r\log T\le x\le c^*t-r\log (t+T)+r\log T+kt^{1/2},$$
even if it means decreasing the positive constant $k$. 

The rest of the proof is as in the homogeneous case.
It follows from (\ref{plan-eq14}) that there exist $x_1>0$ and $x_2>0$, 
both independent of $t \geq h(\tau_0)$ so that 
if we choose $M\ge 1$ large enough then 
\[
Mz(t,c^*t-r\log (t+T)+r\log T+x_1)\ge 2,
\]
and
\[
Mz(t,c^*t-r\log (t+T)+r\log T+x)\le 1/2,~~\hbox{ for all $x>c^*t-r\log (t+T)+r\log T+x_2$.}
\]
Then we set 
\be\label{defbaru}
\bar u(t,x)=\left\{ \begin{array}{cl} 1, & x\le c^*t-r\log (t+T)+r\log T+x_1 \\  \min(1,Mz(t,x)),& x\ge c^*t-r\log (t+T)+r\log T+x_1.\end{array} \right.
\ee
for $t \geq h(\tau_0)$.
Note that $\bar u(t,x)=Mz(t,x)$ for all $x>c^*t-r\log (t+T)+r\log T+x_2$. Moreover, 
$u(0,x) \leq \bar u(h(\tau_0),x)$ for all $x \in \Rm$, even if it means increasing the constant $M$. Therefore, since $\bar u(t,x)$ is a supersolution because of the KPP assumption~(\ref{hypkpp}), the maximum principle implies that
\be\label{ubaru}
u(t,x) \leq \bar u(t + h(\tau_0),x)\ \hbox{ for all }t \geq 0\hbox{ and }x \in \Rm.
\ee
Therefore, for any $\gamma > 0$, we may choose $\bar x$ sufficiently large so that
\[
u(t,x + c^*t - \frac{3}{2\lambda^*} \log(t)) \leq M z(t + h(\tau_0),x + c^*t - \frac{3}{2\lambda^*} \log(t)) < \gamma
\]
holds for all $t > 0$ and $x\ge \bar x$.

Therefore, the proof of the lower bound in Theorem \ref{thm-delay-apr5} is reduced to the proof of Proposition~\ref{thm-plan1}.
The rest of the paper contains the proofs of Propositions~\ref{cor-0524} and \ref{thm-plan1}, as well as that of
Lemma~\ref{prop-0504-2bis}.

%%%%%%%%%%%%%%%%%%%%%%%%%%%%%%%%%%%%%%%%%
%%%%%%%%%%%%%%%%%%%%%%%%%%%%%%%%%%%%%%%%%

\section{Proof of Lemma~\ref{prop-0504-2bis}}

%%%%%%%%%%%%%%%%%%%%%%%%%%%%%%%%%%%%%%%%%

\subsubsection*{The self-adjoint form}

It is useful to write (\ref{2803-8bis}) in a more convenient 
form, to which we can apply the techniques and  ideas of \cite{Norris} where  heat kernel estimates in the whole space are obtained. 

\begin{lem} \label{lem:selfadj} 
Let $\kappa(x)=2\phi_x/\phi$ be defined by (\ref{0523-28}). There is a unique positive, periodic function $\nu(x)$ such that
\begin{equation}\label{0523-26}
\int_0^1\nu(x)dx=1,
\end{equation}
and, for any function $p(x)$,
\begin{equation}
p_{xx}+ \kappa(x) p_x=
\frac{1}{\nu(x)}\pdr{}{x}\left(\nu(x)p_x\right)-\frac{c^*}{\nu(x)}p_x. \label{selfadjId}
\end{equation}
\end{lem}

\noindent{\bf Proof.} The identity (\ref{selfadjId}) means that
\begin{equation}\label{2903-2bis}
\frac{\nu'(x)}{\nu(x)}=\kappa(x)-\farc{\bar b}{\nu(x)},
\end{equation}
with $\bar b = -c^*$, and hence
\begin{equation}\label{0523-20}
\nu_{xx}-\left(\kappa(x)\nu\right)_x=0.
\end{equation}
It is easy to deduce that (\ref{0523-20}) has a positive periodic solution -- this can be seen immediately
since $\tilde\nu(x)\equiv 1$ satisfies the adjoint problem
\[
\tilde\nu_{xx}+\kappa(x)\tilde\nu_x=0,
\]
and by an application of the Krein-Rutman theorem. 

In order to find the constant $\bar b$, observe that the periodic function
\be
\chi(x)=-\farc{1}{\psi(x,\lambda^*)}\frac{d\psi(x,\lambda^*)}{d\lambda} \label{chidef}
\ee
satisfies
\begin{equation}\label{0523-24}
\chi_{xx}+\kappa(x)\chi_x=-\kappa(x)-c^*.
\end{equation}
Indeed, differentiating (\ref{evpperiodic}) in $\lambda$ gives the following equation for $\psi_\lambda=d\psi/d\lambda$:
\[
(\psi_\lambda)_{xx}-2\lambda(\psi_\lambda)_x+\lambda^2\psi_\lambda-2\psi_x+2\lambda\psi+g(x)\psi_\lambda=\gamma'(\lambda)\psi+\gamma\psi_\lambda.
\]
Then, using (\ref{0523-14}) we obtain at $\lambda=\lambda^*$, with $\psi_\lambda^*(x)=\psi_\lambda(x,\lambda^*)$:
\[
(\psi_\lambda^*)_{xx}-2\lambda^*(\psi_\lambda^*)_x+((\lambda^*)^2+g(x))\psi_\lambda^*-2\psi_x^*+2\lambda\psi^*
=c^*\psi+c^*\lambda^*\psi_\lambda^*.
\]
Writing now $\psi_\lambda^*=-\chi(x)\psi(x,\lambda^*)$ and using the definition of $\kappa(x)$
gives (\ref{0523-24}).

Multiplying (\ref{0523-24}) by $\nu(x)$ and integrating over the period gives
\[
\int_0^1\left(\kappa(x)+c^*\right)\nu \,dx=0.
\]
Therefore, we have, since $\nu$  satisfies the normalization (\ref{0523-26}):
\[
-c^*=\int_0^1 \kappa(x)\nu(x)\,dx.
\]
It follows from (\ref{2903-2bis}) that the constant $\bar b$ has to be
$$\bar b=\int_0^1 \kappa(x)\nu(x)\,dx=-c^*.$$
Given that value of $\bar b$, the solution of (\ref{2903-2bis}) exists. \qed
 
\vspace{0.2in}

The periodic function $\chi(x)$ which satisfies (\ref{0523-24}) will be useful later. For this reason, let us remark that there is a unique periodic function $\chi^0(x)$ which satisfies both
$$\chi^0_{xx} + 2 \frac{\phi_x}{\phi} \chi^0_x = -2\frac{\phi_x}{\phi} - c, \quad x \in \Rm$$
and 
\[
\int_0^1 \chi^0(x) \,dx = 0, 
\]
which is obtained by adding a suitable constant to the function $\chi(x)$ defined at (\ref{chidef}). 

%%%%%%%%%%%%%%%%%%%%%%%%%%%%%%%%%%%%%%%%%

\subsubsection*{Proof of Lemma~\ref{prop-0504-2bis}}

We return to (\ref{2803-8bis}) which, by virtue of Lemma \ref{lem:selfadj}, may be written as
\begin{eqnarray}\label{0523-50}
&&p_t=\frac{1}{\nu(x)}(\nu(x)p_x)_x-\farc{c^*}{\nu(x)}p_x, \quad c^*t \le x
\\
&&p(t,c^*t)=0, \quad t > 0\nonumber\\
&&p(0,x)=p_0(x) = u_0(x) e^{\lambda^* x} (\psi(x,\lambda^*))^{-1}, \quad x \geq 0.\nonumber
\end{eqnarray}
Lemma \ref{prop-0504-2bis} is proved using a duality argument, and the main step in the argument is to derive the $L^2$ bound
\begin{equation}\label{2903-48bis}
\left(\int_{c^*t}^\infty p^2(t,x)dx\right)^{1/2}\le 
\frac{C}{t^{3/4}}\int_{0}^\infty xp_0(x)dx, \quad \forall \;t > 0.
\end{equation}

It follows from (\ref{0523-50}) that
\begin{equation}
\farc{1}{2}\farc{d}{dt}\int_{c^*t}^\infty \nu(x)p^2(t,x)dx=
-\int_{c^*t}^\infty\nu(x)p_x^2(t,x)dx. \label{energy1}
\end{equation}
 The right side of (\ref{energy1}) may be bounded from above by using a Nash-type inquality: there is a constant $C$ such that 
\begin{equation}\label{0504-20bis}
\int_0^\infty|\beta(x)|^2dx\le C \left(\int_{0}^\infty \beta_x^2dx\right)^{3/5}
\left(\int_{0}^\infty x\beta (x)dx\right)^{4/5}
\end{equation}
for all functions $\beta \in L^1([0,\infty)) \cap H^1([0,\infty))$ satisfying $\beta(0)=0$ and $\beta(x)\ge 0$ for $x\ge 0$. 
This inequality can be verified in the usual manner: if $\xi(x)$ is an odd extension of
$\beta(x)$ to all of $\Rm$, then 
\begin{equation}\label{0504-22bis}
\int_{-\infty}^\infty|\xi(x)|^2=C\int_{-\infty}^\infty|\hat\xi(k)|^2dk,
\end{equation}
where $\hat\xi(k)$ is the Fourier transform of $\xi(x)$. Note that 
$\hat\xi(0)=0$, and 
\[
\left|\frac{d}{dk}\hat\xi(k)\right|\le C\int_0^\infty x\beta(x)dx,
\]
whence $|\hat\xi(k)|\le C|k|\|x\beta\|_1$. It follows from (\ref{0504-22bis}) that for any $R>0$ we have
$$\int_{-\infty}^\infty|\xi(x)|^2dx\le C\int_{|k|\le R}|\hat\xi(k)|^2dk+
C\int_{|k|\ge R}\frac{|k|^2}{R^2}|\hat\xi(k)|^2dk
\le CR^3\|x\beta\|_1^2+\frac{C}{R^2}\|\beta_x\|_2^2.$$
Choosing $R=\left(\|\beta_x\|_2^2 \;/\; \|x\beta\|_1^2\right)^{1/5}$ gives (\ref{0504-20bis}).   

Going back to (\ref{energy1}), since $\nu(x)^{-1} > 0$ is bounded, we conclude that
\begin{equation}
\farc{1}{2}\farc{d}{dt}\int_{c^*t}^\infty \nu(x)p^2(t,x)dx \le - C \left(\int_{c^*t}^\infty (p(t,x))^2 \,dx \right)^{5/3} 
\left(\int_{c^*t}^\infty (x-c^*t) p(t,x) \,dx\right)^{-4/3}. \label{nash1}
\end{equation}
Next, we work toward an estimate of the right side of (\ref{nash1}). 
Let us multiply (\ref{0523-50}) by a function $\nu(x)f(t,x)$, with $f(t,c^*t)=0$ and integrate:
$$\begin{array}{rcl}
\displaystyle\frac{d}{dt}\int_{c^*t}^\infty \nu(x)f(t,x)p(t,x)\,dx & = & \displaystyle\int_{c^*t}^\infty \nu(x)f_t(t,x)p(t,x)\,dx
-\int_{c^*t}^\infty \nu(x)f_x(t,x)p_x(t,x)\,dx\vspace{3pt}\\
& & \displaystyle- c^*\int_{c^*t}^\infty fp_x \,dx.\end{array}$$
We will choose $f$ to be a solution of the backward equation, defined by the next lemma.

\begin{lem} \label{lem:fexist}
There is a function $f(t,x)$ and a constant $m > 0$ such that $f_t<0$,
\be\label{0524-18}\left\{\begin{array}{l}
f_t+\displaystyle\farc{1}{\nu(x)}\left(\nu(x)f_x\right)_x+\farc{c^*}{\nu(x)}f_x=0, \quad x > c^*t, \quad t \in \Rm, \\
f(t,c^* t) = 0, \quad t \in \Rm
\end{array}\right.
\ee
and 
\begin{equation}\label{2903-24bis}
m(x-c^*t)<f(t,x)<m^{-1}(x-c^*t),~~\hbox{for all $x>c^*t$, $t \in \Rm$.}
\end{equation}
\end{lem}

Let us postpone the proof of this lemma for the moment, and use it to finish the proof of Lemma~\ref{prop-0504-2bis}.
Given the function $f(t,x)$ described in Lemma \ref{lem:fexist}, observe that the integral
$$I(t)=\int_{c^*t}^\infty \nu(x)f(t,x)p(t,x)dx$$
is preserved: $I(t)=I(0)$ for all $t \geq 0$. Moreover, (\ref{2903-24bis}) implies that
\[
\left(~\int_{c^*t}^\infty (x - c^*t) p(t,x) \,dx\right)^{-4/5} \geq C\left(~\int_{c^*t}^\infty \nu(x)f(t,x)p(t,x)\,dx\right)^{-4/5} = C(I(0))^{-4/5},
\]
for all $t > 0$. Therefore, if
\[
I_2(t)=\int_{c^*t}^\infty \nu(x) p^2(t,x)dx,
\]
then from (\ref{nash1}) we conclude
\[
\frac{dI_2(t)}{dt} \le -C\frac{(I_2(t))^{5/3}}{(I(0))^{4/3}}.
\]
It follows that $(I_2(t))^{-2/3} \ge Ct (I(0))^{-4/3}$ for all $t > 0$, which implies the $L^2$ bound (\ref{2903-48bis}). 

The standard duality argument can be now applied. If $S_t$ is the solution operator mapping $p_0(\cdot)$ to $p(t,\cdot)$, then the  adjoint operator $S_t^*$ is of the same form as $S_t$ except for 
$c^*$ replaced by $(-c^*)$ and changing the direction of time.
Hence,   the $L^1\to L^2$ bound (\ref{2903-48bis}) for $S_t$ implies also the dual $L^2\to L^\infty$ bound:
\[
|p(t,x)|\le \farc{C(x-c^*t)}{t^{3/4}}\|p_0\|_{L^2}, \quad x > c^*t, \;\; t > 0.
\]
Finally, writing $S_t=S_{t/2}\circ S_{t/2}$ we obtain the conclusion of Lemma~\ref{prop-0504-2bis}.~$\Box$

\vspace{0.2in}

\noindent{\bf Proof of Lemma \ref{lem:fexist}.} Observe that (\ref{0524-18}) has a solution of the form $Y(t,x) = (x - c^*t) + y(x)$, where $y(x)$ is periodic and satisfies
\[
-c^*\nu(x) +\left(\nu(x)(1+y_x)\right)_x+c^*(1+y_x)=0,
\]
or
\begin{equation}\label{2903-20bis}
\left(\nu(x)y_x\right)_x+c^*y_x=c^*(\nu(x)-1)-\nu'(x).
\end{equation}
Equation (\ref{2903-20bis}) has a periodic solution because the integral of the right side over the period vanishes, 
because of (\ref{0523-26}). By subtracting a constant from $y$, we may assume $Y(t,c^*t) \leq 0$. Although $Y(t,x)$ grows linearly in $(x - c^*t)$ and is a solution of (\ref{0524-18}) 
for all $t \in \Rm$ and $x \in \Rm$, it may  not satisfy the desired Dirichlet boundary condition at $x = c^*t$. 
On the other hand, if $\beta(t)$ is the largest zero of $Y$ then 
\begin{equation}\label{2903-36bis}
|\beta(t)-c^*t|\le M,
\end{equation}
with a constant $M$ that does not depend on $t$.

A function $f(t,x)$ having the desired properties may be constructed as the limit of the sequence of functions $f^{(n)}(t,x)$ which satisfy
\begin{eqnarray}
&&f^{(n)}_t+\farc{1}{\nu(x)}\left(\nu(x)f^{(n)}_x\right)_x+\farc{c^*}{\nu(x)}f^{(n)}_x=0, \quad x > c^*t, \quad t \leq n \no \\
&& f^{(n)}(t,c^*t) = 0 \quad t \leq n, \no \\
&& f^{(n)}(n,x) = \max(0,Y(n,x)), \quad x \geq c^*n. \no
\end{eqnarray}
It follows from the maximum principle and (\ref{2903-36bis}) that
there exists a constant $C$, independent of $n$, such that
\begin{equation}\label{2903-42bis}
 Y(t,x) - C \le f^{(n)}(t,x) \le Y(t,x) + C, \quad \forall \; x \geq c^*t, \; t \leq n.
\end{equation}
Using (\ref{2903-42bis}), we can find positive constants $L$, $M$, $m$, independent of $n$, so that 
$$\hbox{$f^{(n)}(t,ct+L)>M_1$, for all $t\le n$},$$
and, in addition, $m(x-c^*t)<f^{(n)}(t,x)<m^{-1}(x-c^*t)$ holds for $x>c^*t+L$ and $t < n/2$. Then the strong maximum principle and parabolic regularity
imply that $f^{(n)}_x(t,c^*t)>c_0$ for all $t<n/2$, for some positive constant $c_0$ that does not depend on $n$ or $t$. 
By parabolic regularity, we may then extract a subsequence converging to a limit $f(t,x)$ satisfying (\ref{0524-18}), (\ref{2903-24bis}) 
and the boundary condition $f(t,c^*t) = 0$ for all $t \in \Rm$.  Note that  $f^{(n)}_t\le 0$  -- this follows from the maximum 
principle since $f^{(n)}(t,x)\ge 0$ and $f^{(n)}(t,x)\ge Y(t,x)$ for all $t\le n$, and $x\ge c^*t$.  
It follows that in the limit we also have $f_t(t,x)\le 0$. \qed

%%%%%%%%%%%%%%%%%%%%%%%%%%%%%%%%%%%%%%%%%
%%%%%%%%%%%%%%%%%%%%%%%%%%%%%%%%%%%%%%%%%

\section{The proof of Proposition \ref{cor-0524}}

Proposition \ref{cor-0524} is based on the following key estimate, which is proved in Section \ref{sec:prop-0523}.

\begin{prop}\label{prop-0523}
There exist a time $T_0>0$ and constants $c_0>0$, $\beta > 0$, and $N>0$
that depend only on the initial data so that for any $t>T_0$ there exists a set $I_t \subset [c^*t + N^{-1} \sqrt{t} , c^*t+N\sqrt{t}]$ such that $|I_t| \geq \beta \sqrt{t}$ and
\begin{equation}\label{0523-52}
p(t,x)\ge\frac{c_0}{t}.
\end{equation}
holds for all $x \in I_t$.
\end{prop}

We also make use of an estimate for the heat kernel associated with the equation
\begin{equation}
\rho_t = \frac{1}{\nu(x)}(\nu(x) \rho_x)_x - \frac{c^*}{\nu(x)} \rho_x. \label{feqnline}
\end{equation}
For $R > 0$ and $\xi \in \Rm$ fixed, let $\bar \Gamma(t,x,s,y) = \bar \Gamma(t,x,s,y;R,\xi)$ denote the heat kernel for (\ref{feqnline}) in the tilted cylinder
\[
T(\xi,R,s) = \left \{ (t,x) \in \Rm^2\;:\;\; |x - \xi - c^*t| < R , \quad t \geq s \right\}
\]
with the Dirichlet boundary conditions on the lateral boundary of the cylinder. 
That is, if $s \in \Rm$ and $|y - \xi - cs| < R$, $\bar \Gamma(t,x,s,y)$ satisfies (\ref{feqnline}) for $(t,x) \in T(\xi,R,s)$, with 
the boundary condition $\bar \Gamma(t,x,s,y) = 0$ if $|x - \xi - c^*t| = R$, and the initial condition 
\[
\lim_{t \searrow s} \bar \Gamma(t,x,s,y) = (\nu(y))^{-1} \delta_y(x).
\]
The following lemma gives a lower bound on $\bar \Gamma(t,x,s,y)$, provided that $x$ and $y$ are sufficiently far from the boundary of $T(\xi,R,s)$.

\begin{lem} \label{lem:heatlower}
For all $\delta \in (0,1)$, there are some constants $\alpha > 0$ and $K>0$ such that
$$\bar \Gamma(t,x,s,y - c^*(t-s);R,\xi)  \geq \frac{\alpha}{2K(t-s)^{1/2}}  e^{-\frac{K|y-x|^2}{(t-s)}}$$
holds if $R > 0$, $t \in (s, s + R^2]$, and $x,y \in (c^*t + \xi - \delta R, c^* t + \xi + \delta R)$.
\end{lem}

\noindent{\bf Proof.} Let $\Gamma(t,x,s,y)$ denote the free-space heat kernel associated with the equation 
$$\nu(x) \rho_t = (\nu(x) \rho_x)_x - c^* \rho_x.$$
That is, for each $(s,y) \in \Rm^2$, $\Gamma(\cdot, \cdot,s,y)$ is a solution of (\ref{feqnline}) for $x \in \Rm$ and $t > s$, with
\[
\lim_{t \searrow s} \Gamma(t,x,s,y) = (\nu(y))^{-1}\delta_{y}(x).
\]
If $\rho$ is continuous, bounded and satisfies (\ref{feqnline}) for $x \in \Rm$ and $t > s$, then
\[
\rho(t,x) = \int_{\Rm} \Gamma(t,x,s,y) \rho(s,y) \nu(y) \,dy
\]
for $t \geq s$. For $\Gamma(t,x,s,y)$ we have the following estimates of Norris \cite{Norris}, Theorem 1.1: 
there is a constant $K > 0$ such that
\begin{equation}
\farc{e^{-K|x-y|^2/(t-s)}}{K |t-s|^{1/2}}   \leq \Gamma(t,x,s,y - c^* (t-s)) \leq \farc{Ke^{-|x-y|^2/K(t-s)}}{ |t-s|^{1/2}} . 
\label{norrisbound}
\end{equation}
holds for all $x,z \in \Rm$, $t > s$.  Obviously, (\ref{norrisbound}) implies the upper bound
\[
\bar \Gamma(t,x,s,y - c^*(t-s);R,\xi) \leq \Gamma(t,x,s,y - c^*(t-s)) \leq K |t-s|^{-1/2} e^{-|x-y|^2/K(t-s)}.
\]

The proof of Lemma \ref{lem:heatlower} mimics the analysis of Fabes and Stroock \cite{FS} (see the proof Lemma 5.1, therein). 
It suffices to assume $s = 0$ and $\xi = 0$. The first step is to derive the identity
\begin{equation}
\bar \Gamma(t,x,0,y) =  \Gamma(t,x,0,y) -  \int_0^t \left( \Gamma(t,x,r,c^*r + R)h^+(r) +  \Gamma(t,x,r,c^*r-R)h^-(r) \right) \,dr \label{greenident3}
\end{equation}
where $h^{\pm}(r) \geq 0$ depends on $y$ and $R$, but 
\[
\int_0^t (h^+(r) + h^-(r)) \,dr \leq 1
\]
always holds. This is analogous to a statement on p. 335 of \cite{FS}.  To see where (\ref{greenident3}) comes from, 
suppose $\rho$ satisfies (\ref{feqnline}) for $(t,x) \in T_R = \{ (t,x) \;|\; t \geq 0,\;\; |x - c^*t| < R \}$. 
Choose a test function $\varphi(r,z)$, and integrate over $r \in [t_1,t_2]$, $z \in D_r = [c^*r-R, c^*r + R]$:
$$\begin{array}{rcl}
0 & = & \displaystyle\int_{t_1}^{t_2} \!\! \int_{D_r}  \left(\nu(z) \rho_r - (\nu(z) \rho_z)_z + c^* \rho_z \right) \varphi(r,z) \,dz \,dr \no \\
& = & \displaystyle- \int_{t_1}^{t_2} \!\! \int_{D_r} \rho(r,z) \left(\nu \varphi_r + (\nu \varphi_z)_z - c^* \varphi_z \right) \,dz \,dr 
 +  \int_{t_1}^{t_2} \!\! \int_{D_r}  (\nu \rho \varphi_z)_z -(\nu \rho_z \varphi)_z  \,dr \,dz \no \\
&  & \displaystyle+  \int_{t_1}^{t_2} \!\! \int_{D_r} (\nu \rho \varphi)_r + (c^* \rho \varphi)_z \,dr \,dz.
\end{array}$$
So, if $\varphi$ satisfies the adjoint equation $\nu \varphi_r + (\nu \varphi_z)_z - c^* \varphi_z = 0$ 
and if $\rho$ vanishes on $\partial D_r$ for each $r$, we obtain
\begin{eqnarray}
 \int_{D_{t_2}} \nu \rho \varphi   \,dz - \int_{D_{t_1}} \nu \rho \varphi   \,dz  = \int_{t_1}^{t_2} \!\! \nu \rho_z \varphi \Big|_{c^*r-R}^{c^*r+R}   \,dr.  \label{greenident}
\end{eqnarray}
Now, let $t > 0$ and $x \in D_t = (c^*t - R,c^*t + R)$ be fixed. For $y \in (-R,R)$ fixed, let $\rho(r,z) = \bar \Gamma(r,z,0,y)$ 
and $\varphi(r,z) = \Gamma(t,x,r,z)$. The function $\rho(r,z)$ satisfies (\ref{feqnline}) in $T(\xi,R,s)$ with $\rho(r,z) = 0$ for $z \in D_r$. 
The function $\varphi(r,z)$ is a solution of the adjoint equation $\nu \varphi_r + (\nu \varphi_z)_z - c \varphi_z = 0$ for $r \in (0,t)$. 
Therefore, (\ref{greenident}) holds. By letting $t_2 \to t$ and $t_1 \to 0$, we obtain the following identity relating $\bar \Gamma$ 
to the free-space heat kernel $\Gamma$:
$$\bar \Gamma(t,x,0,y) =  \Gamma(t,x,0,y) +  \int_0^t \nu(z) \bar \Gamma_z(r,z,0,y) \Gamma(t,x,r,z) \Big|_{z=c^*r-R}^{z=c^* r+R}   \,dr.$$
Here we have used the fact that
\[
\lim_{t_2 \nearrow t} \int \Gamma(t,x,t_2,z)f(z)\nu(z) \,dz = f(x) \quad \quad \text{and} 
\quad \quad \lim_{t_1 \searrow 0} \int \bar \Gamma(t_1,z,0,y)f(z)\nu(z) \,dz = f(y)
\]
for any continuous $f$. Note that $\bar \Gamma_z(r,c^*r+R,0,y) \leq 0$ and  $\bar \Gamma_z(r,c^*r-R,0,y) \geq 0$. 
If we had chosen $\varphi \equiv 1$, instead, we would have obtained
$$\begin{array}{l}
\displaystyle\int_{D_t} \nu(x) \bar \Gamma(t,x,0,y) \,dx\\
\qquad\qquad = \displaystyle1 +  \int_0^t \!\! \nu(z) \bar \Gamma_z(r,z,0,y)  \Big|_{z=c^*r-R}^{z=c^*r+R}   \,dr  \no \\
\qquad\qquad=\displaystyle1 -  \int_0^t \!\! \left(\nu(c^*r + R)  |\bar \Gamma_z(r,c^*r + R,0,y)| + \nu(c^*r-R) |\bar \Gamma_z(r,c^*r -R ,0,y)| \right)   \,dr. \no
\end{array}$$
Since the left side is non-negative, this implies 
\[
 \int_0^t \!\! \left(\nu(c^*r + R)  |\bar \Gamma_z(r,c^*r + R,0,y)| + \nu(c^*r-R) |\bar \Gamma_z(r,c^*r -R,0,y)| \right)   \,dr \leq 1.
\]
Thus, we have shown (\ref{greenident3}).

By combining (\ref{greenident3}) with the estimate (\ref{norrisbound}) for $\Gamma$, we obtain a lower bound on $\bar \Gamma$:
\begin{equation}
\bar \Gamma(t,x,0,y - c^*t) \geq \frac{e^{-K|y-x|^2/t}}{K t^{1/2}}   -  
K \sup_{0 \leq \tau \leq t} \farc{e^{-R^2(1 - \delta)^2/(K \tau)}}{\tau^{1/2}} ,\label{greenident4}
\end{equation}
for all $x \in [-\delta R,\delta R]$, $y \in [-R,R]$, $t > 0$. Observe that the unique maximum of the function
\[
\beta(\tau)=\farc{ e^{-R^2(1 - \delta)^2/(K \tau)}}{\tau^{1/2}} , \quad \tau > 0,
\]
occurs at the point $\tau^* = 2 R^2(1 - \delta)^2/K$. So, if $\varepsilon^2 < 2(1 - \delta)^2/K$ and $t \leq \varepsilon^2 R^2$, we have $t \leq \tau^*$. In this case, (\ref{greenident4}) gives us the bound
$$\begin{array}{rcl}
\bar \Gamma(t,x,0,y - c^*t) & \geq & \displaystyle\farc{e^{-K|y-x|^2/t}}{K t^{1/2}}   -  K \sup_{0 \leq \tau \leq t} \farc{ e^{-R^2(1 - \delta)^2/(K \tau)}}{\tau^{1/2}} \no \\
& = & \displaystyle\farc{e^{-K|y-x|^2/t}}{K t^{1/2}} -  K   \farc{ e^{-R^2(1 - \delta)^2/(K t)}}{t^{1/2}}  \\
& = & \displaystyle\farc{e^{-K|y-x|^2/t}}{K t^{1/2}} \left(  1  -  K^2   e^{- \frac{R^2(1 - \delta)^2}{K t} + \frac{K|x-y|^2}{t}}  \right).
\end{array}$$
If $x \in [-\delta R,\delta R]$ and $|x - y| \leq \varepsilon R$ also hold, 
and $\varepsilon^2 < (1 - \delta)^2/(2K^2)$ is small enough we have
$$\begin{array}{rcl}
1  -  K^2   e^{- \frac{R^2(1 - \delta)^2}{K t} + \frac{K|x-y|^2}{t}} & \geq & 
 1  -  K^2   e^{-t^{-1}( \frac{R^2(1 - \delta)^2}{K } - K\varepsilon^2 R^2)}
 \geq  1  -  K^2   e^{-t^{-1} \frac{R^2(1 - \delta)^2}{2K }} \\
& \geq & 1  -  K^2   e^{-\varepsilon^{-2}( \frac{(1 - \delta)^2}{2K })} > 1/2.
\end{array}$$
This implies that for any $\delta \in (0,1)$ and $R > 0$,
$$\bar \Gamma(t,x,0,y - c^*t) \geq \frac{1}{2Kt^{1/2}}  e^{-\frac{K|y-x|^2}{t}}$$
if $x \in [-\delta R,\delta R]$ and $|x - y| \leq \varepsilon R$, $t \leq \varepsilon^2 R^2$, and $\varepsilon$ is sufficiently small, 
depending only on $\delta$ and~$K$.  A chaining argument, as in \cite{FS},
now shows that for any $\delta \in (0,1)$, 
there must be a constant $\alpha$, depending only on $\delta$ and $K$, such that
$$\bar \Gamma(t,x,0,y - c^*t)  \geq \frac{\alpha}{2Kt^{1/2}}  e^{-\frac{K|y-x|^2}{t}}$$
holds if $x,y \in [-\delta R,\delta R]$, $t \leq R^2$ (i.e. rather than just $t \leq \varepsilon^2 R^2$). 
Although $\bar \Gamma$ depends on $R$, $\alpha$ and $K$ are independent of $R$. This finishes the proof of Lemma~\ref{lem:heatlower}.
\qed

%%%%%%%%%%%%%%%%%%%%%%%%%%%%%%%%%%%%%%%%%

\subsubsection*{End of the proof of Proposition~\ref{cor-0524}}

We may now finish the proof of Proposition \ref{cor-0524}.  By Proposition \ref{prop-0523} we have 
\begin{equation}\label{0729-02}
p(s,x)\geq \dfrac{c_0}{ s}
\end{equation}
for all $s \geq T_0$ and $x \in I_s$, where $I_s \subset [c^*s+N^{-1}\sqrt s,c^*s+N\sqrt s]$ and $|I_s| \geq \beta \sqrt{s}$. 
We apply the lower bound on the heat-kernel   in Lemma \ref{lem:heatlower}. 
Let $s \geq T_0$, $R = \sqrt{s}(N^{-1} + N)/2$, $\xi = c^*s + R$, and 
$\bar \Gamma = \bar \Gamma(t,x,s,y;R,\xi)$ be the heat kernel in the tilted cylinder $T(\xi,R,s)$ with 
Dirichlet boundary conditions. For $x \in [c^*t ,c^*t + 2R]$, $t > s$, we have
\begin{eqnarray}
p(t,x) \geq \int_{cs}^{cs + 2R} \bar \Gamma(t,x,s,y)p(s,y) \nu(y)\,dy. \label{pkernellower}
\end{eqnarray}
Set 
\[
\delta =  \frac{N-N^{-1}}{N+N^{-1} } \in (0,1),
\]
and $t = s + R^2$. Observe that 
\[
I_s \subset [c^*s+N^{-1}\sqrt s,c^*s+N\sqrt s] = [c^*s + (1 - \delta)R, c^*s + (1 + \delta) R].
\]
By Lemma \ref{lem:heatlower}, we have
\[
\bar \Gamma(t,x,s,y) \geq \frac{\alpha}{2 (t -s)^{1/2}} e^{- \frac{K|x - y|^2}{t - s}} =  \frac{\alpha}{2 KR} e^{- \frac{K|x - y|^2}{R^2}}
\]
for all 
\[
x \in [c^*t + (1 - \delta)R ,c^*t + (1 + \delta)R] = [c^*t + N^{-1} \sqrt{s}, c^*t + N\sqrt{s}],
\]
and 
\[
y \in [c^*s + (1 - \delta)R ,c^*s + (1 + \delta)R] = [c^*s + N^{-1} \sqrt{s}, c^*s + N\sqrt{s}].
\]
Therefore, by combining (\ref{0729-02}) and (\ref{pkernellower}) we obtain
$$\begin{array}{rcl}
p(t,x) & \geq & \displaystyle\int_{I_s} \bar \Gamma(t,x,s,y)p(s,y) \nu(y)\,dy\\
&  \geq & \displaystyle|I_s|   \min_{y \in I_s} \bar \Gamma(t,x,s,y)p(s,y) \nu(y)  
 \geq   |I_s|   \frac{C}{\sqrt{s}} \min_{y \in I_s} p(s,y)  \geq \frac{C}{s},\end{array}$$
for all $x \in [c^*t + (1 - \delta)R ,c^*t + (1 + \delta)R]$. Since $R = \sqrt{s}(N^{-1} + N)/2$ and 
$t = s + R^2$ we have shown that for $\sigma = 1 + (N^{-1} + N)^2/4$, there is a constant $C > 0$ such that
$$p(t,c^* t + \sigma \sqrt{t}) \geq  \frac{C}{ s} = \frac{C \sigma}{ t}$$
holds for all $t \geq \sigma T_0$. Therefore, the last remaining ingredient in the proof of the lower bound in Theorem~\ref{thm-delay-apr5} 
is the proof of Proposition~\ref{prop-0523}.~$\Box$

%%%%%%%%%%%%%%%%%%%%%%%%%%%%%%%%%%%%%%%%%
%%%%%%%%%%%%%%%%%%%%%%%%%%%%%%%%%%%%%%%%%

\section{The proof of Proposition~\ref{prop-0523}} \label{sec:prop-0523}

%%%%%%%%%%%%%%%%%%%%%%%%%%%%%%%%%%%%%%%%%

\subsection{The homogeneous case}

Since the proof of Proposition~\ref{prop-0523} is rather long we first present it in the simplest case $\nu(x)\equiv 1$,
$c^*=0$. In that case (\ref{0523-52}) (and (\ref{pexplheat})) can be proved simply by examining the explicit formula for the solution
to the heat equation on the half line, as shown in \cite{HNRR1}. However, as such formulas are not available in the non-uniform case, we will
present an alternative (and much longer!)  proof using  the energy  method that we will adapt to the periodic case. 
The key step is the following lemma.

\begin{lem}\label{lm-0805bis} Let $\nu(x)\equiv 1$ and $c^*=0$, and let $p(t,x)$ solve (\ref{0523-50}) with $p_0(x)$ being compactly supported on $[0,\infty)$. There exists $C>0$ so that for any $\alpha>0$ we have
\begin{equation}\label{080502bis}
\left(\int_0^\infty\farc{e^{2\alpha x}-e^{-2\alpha x}}{2\alpha x}p^2(t,x)dx\right)^{1/2}\le 
\frac{Ce^{\alpha^2t}}{t^{3/4}}
\int_0^\infty\farc{e^{\alpha x}-e^{-\alpha x}}{\alpha}p_0(x)dx.
\end{equation}
\end{lem}

Let us first show how (\ref{0523-52}) follows from (\ref{080502bis}). We will take $\alpha=1/\sqrt{t}$ in (\ref{080502bis}).
Then, if $T_0$ is sufficiently large, and $t>T_0$, for any $x\in\hbox{supp}~p_0$ we have
\[
\farc{e^{\alpha x}-e^{-\alpha x}}{\alpha}\le 4x.
\]
Moreover, the integral
\[
I(t)=\int_0^\infty xp(t,x)dx 
\]
is conserved: $I(t)=I(0)$. We conclude that for all $t>T_0$ we have
\begin{equation}\label{080506bis}
\left(\int_0^\infty\farc{e^{2 x/\sqrt{t}}-e^{-2x/\sqrt{t}}}{2x/\sqrt{t}}p^2(t,x)dx\right)^{1/2}\le 
\frac{C}{t^{3/4}}
\int_0^\infty xp_0(x)dx,
\end{equation}
or
\begin{equation}\label{080508bis}
\left(\int_0^\infty\farc{e^{2 x/\sqrt{t}}-e^{-2x/\sqrt{t}}}{x}p^2(t,x)dx\right)^{1/2}\le 
\frac{C}{t}
\int_0^\infty xp_0(x)dx.
\end{equation}
Let us now take $N>1$ sufficiently large (but independent of $t$), then for $x>N\sqrt{t}$ we have 
$e^{2x/\sqrt{t}}>2e^{-2x/\sqrt{t}}$, thus (\ref{080508bis}) implies
$$\left(\int_{N\sqrt{t}}^\infty\farc{e^{2 x/\sqrt{t}}}{x}p^2(t,x)dx\right)^{1/2}\le 
\frac{C}{t}
\int_0^\infty xp_0(x)dx.$$
Moreover, we have
$$\begin{array}{rcl}
\displaystyle\int_{N\sqrt{t}}^\infty xp(t,x)dx & \le & \displaystyle\int_{N\sqrt{t}}^\infty \frac{e^{x/\sqrt{t}}}{\sqrt{x}}p(t,x)e^{-x/\sqrt{t}}x^{3/2}dx\\
& \le & \displaystyle\left(\int_{N\sqrt{t}}^\infty \frac{e^{2x/\sqrt{t}}}{{x}}p^2(t,x)dx\right)^{1/2}
\left(\int_{N\sqrt{t}}^\infty e^{-2x/\sqrt{t}}x^{3}dx\right)^{1/2}\\
& \le & \displaystyle{C}\left(\int_0^\infty xp_0(x)dx\right)\left(\int_N^\infty y^3e^{-y}dy\right)^{1/2}\\
& \le & \displaystyle\farc{1}{4}\int_0^\infty xp_0(x)dx=\farc{I(0)}{4},
\end{array}$$
as long as $N>N_0$ is large enough (but independent of $t$). As $I(t)=I(0)$, it follows that 
\[
\int_0^{N\sqrt{t}}xp(t,x)dx \ge \frac{3I_0}{4}.
\]
From Lemma \ref{prop-0504-2bis} we know that
\[
\int_0^{N^{-1}\sqrt{t}}xp(t,x)dx \leq C N^{-3} I_0.
\]
Therefore, by taking $N$ larger, if necessary, we have
\[
\int_{N^{-1} \sqrt{t}}^{N\sqrt{t}}xp(t,x)dx \ge \frac{I_0}{2}.
\]
For $c_0 > 0$ to be chosen, let $H_t^\pm$ be the sets $H_t^+ = \{ x \in [N^{-1} \sqrt{t}, N \sqrt{t}] \;|\; p(t,x) \geq c_0/t \}$, 
and $H_t^- = \{ x \in [N^{-1} \sqrt{t}, N \sqrt{t}] \;|\; p(t,x) < c_0/t \}$. We have
\[
\frac{I_0}{2} \leq \int_{H_t^+} xp(t,x)dx + \int_{H_t^-} xp(t,x)dx \leq \int_{H_t^+} xp(t,x)dx  + \frac{c_0}{2} N^2.
\]
so that by choosing $c_0 \leq I_0/(2N^2)$, we have
\[
\frac{I_0}{4} \leq \int_{H_t^+} xp(t,x)dx.
\]
Now, apply Lemma \ref{prop-0504-2bis} again:
\[
\frac{I_0}{4}  \leq \int_{H_t^+} xp(t,x)dx \leq \frac{C I_0}{t^{3/2}} \int_{H_t^+} x^2 dx \leq \frac{C I_0}{t^{3/2}} |H_t^+| N^2 t.
\]
It follows that $|H_t^+| \geq \sqrt{t}/(4 N^2 C)$. Except
for the proof of Lemma~\ref{lm-0805bis}, this proves Proposition~\ref{prop-0523} in the homogeneous case. 

%%%%%%%%%%%%%%%%%%%%%%%%%%%%%%%%%%%%%%%%%

\subsubsection*{\bf Proof of Lemma~\ref{lm-0805bis}} 

In the homogeneous case  (\ref{0523-50}) is simply
\[
p_t = p_{xx},~~p(t,0)=0,
\]
since we assume that $\nu(x) \equiv 1$ and $c^* = 0$. There are at least three ways to prove Lemma~\ref{lm-0805bis} in 
this situation: first, one can use the explicit formula for $p(t,x)$. Second, one can use the fact that 
$q(t,x)=p(t,x)/x$ solves 
\begin{equation}\label{270404bis}
q_t=q_{xx}+\farc{2}{x}q_x.
\end{equation}
Hence, if we set $\tilde q(t,z)=q(t,|z|)$, with $z\in\Rm^3$, then we get the heat equation in $\Rm^3$ for $\tilde q$:
\begin{equation}\label{270406bis}
\tilde q_t=\Delta_z \tilde q,~~z\in\Rm^3.
\end{equation}
Then one could apply the usual Nash inequality to the function
\[
\varphi(t,z)=e^{\bar\alpha\cdot z}\tilde q(t,z),~~~z\in\Rm^3,
\]
and prove Lemma~\ref{lm-0805bis} in this way.
Neither of these methods would generalize to the periodic case, hence we develop a third, longer but generalizable proof. Motivated by the above, let us define the exponential moments 
\begin{eqnarray}\label{270410bis}
I_\alpha(t)= \frac{1}{2 \pi} \int_{\Rm^3} \varphi(t,z) \,dz=\int_0^\infty\farc{e^{\alpha x}-e^{-\alpha x}}{\alpha}xq(t,x)dx,
\end{eqnarray}
$$V_\alpha(t) = \frac{1}{2\pi} \int_{\Rm^3} (\varphi(t,z))^2 \,dz=\int_0^\infty\farc{e^{2\alpha x}-e^{-2\alpha x}}{2\alpha}xq^2(t,x)dx,$$
and 
$$D_\alpha(t)=\frac{1}{2\pi}\int_{\Rm^3}  |\nabla\varphi|^2 \,dz
=\int_0^\infty\farc{e^{2\alpha x}-e^{-2\alpha x}}{2\alpha}(q_x^2-\alpha^2q^2)xdx.$$
The Nash inequality in $\Rm^3$ (e.g. \cite{Stroock}, Lemma I.1.1.)  gives the following lemma.

\begin{lem} \label{lem:NashR3}
There is a constant $C > 0$ such that for any function $w(x):[0,\infty) \to \Rm$ which is smooth, bounded and compactly supported we have
\[
(\hat V_\alpha)^{5/3} \leq C (\hat I_\alpha)^{4/3} \hat D_\alpha,
\]
for all $\alpha > 0$, where
\[
\hat V_\alpha = \int_0^\infty \frac{e^{2\alpha x} - e^{-2\alpha x}}{2\alpha} x w^2(x) \,dx,
~~~~
\hat I_\alpha = \int_0^\infty \frac{e^{\alpha x} - e^{-\alpha x}}{\alpha} x |w(x)| \,dx,
\]
and
\[
\hat D_\alpha = \int_0^\infty \frac{e^{2\alpha x} - e^{-2\alpha x}}{2\alpha} x(w_x^2 - \alpha^2 w^2) \,dx.
\]
\end{lem}

Using (\ref{270406bis}), it is easy to check that $I_\alpha(t) = e^{\alpha^2 t} I_\alpha(0)$ and 
\begin{equation}
V_\alpha'(t) = 2 \alpha^2 V_\alpha(t) - 2D_\alpha(t).\label{Valphadot}
\end{equation}
Lemma~\ref{lem:NashR3} applied to (\ref{Valphadot}) results in the bound
\[
V_\alpha'(t)\le 2\alpha^2 V_\alpha(t)-\frac{C[V_\alpha(t)]^{5/3}}{[I_\alpha(t)]^{4/3}}.
\]
If $V_\alpha(t)=e^{2\alpha^2t}Z(t)$, then
\[
Z'(t)\le -\frac{CZ_\alpha(t)^{5/3}}{(I_\alpha(0))^{4/3}}.
\]
It follows that $Z(t)\le {(I_\alpha(0))^{2}}/{t^{3/2}},$ hence
$$\left(\int_0^\infty\farc{e^{2\alpha x}-e^{-2\alpha x}}{2\alpha}xq^2(t,x)dx\right)^{1/2}\le \frac{Ce^{\alpha^2t}}{t^{3/4}}
\int_0^\infty\farc{e^{\alpha x}-e^{-\alpha x}}{\alpha}xq(0,x)dx.$$
In terms of $p(t,x)$ this is:
$$\left(\int_0^\infty\farc{e^{2\alpha x}-e^{-2\alpha x}}{2\alpha x}p^2(t,x)dx\right)^{1/2}\le 
\frac{Ce^{\alpha^2t}}{t^{3/4}}
\int_0^\infty\farc{e^{\alpha x}-e^{-\alpha x}}{\alpha}p_0(x)dx.$$
This completes the proof of Lemma~\ref{lm-0805bis}. $\Box$

%%%%%%%%%%%%%%%%%%%%%%%%%%%%%%%%%%%%%%%%%

\subsection{The general case}

We now adapt the preceding proof to the general problem
\begin{equation}\label{090502}
p_t=\frac{1}{\nu(x)}\pdr{}{x}\left(\nu(x)p_x\right)-\farc{c^*}{\nu(x)}p_x,
\end{equation}
with the Dirichlet boundary condition $p(t,c^*t)=0$.  The next lemma gives the analog of the function $x$
in the periodic case.

\begin{lem} \label{lem:zetabounds}
There is a function $\zeta(t,x)$ and a constant $m > 0$ such that
\begin{eqnarray}
&&\nu(x) \zeta_t= \left(\nu(x)\zeta_x\right)_x - c^*\zeta_x,\quad x>c^*t,~~t \in \Rm,\label{090504} \\
&&\zeta(t,c^*t)=0, \quad t \in \Rm\label{zetabc}
\end{eqnarray}
and
$$m (x - c^*t) \leq \zeta(t,x) \leq m^{-1} (x - c^*t), \quad \forall \; x \geq c^*t, \quad t \in \Rm.$$
\end{lem}

In analogy to the uniform case, define
\begin{equation}\label{march204}
\hbox{$q(t,x)={p(t,x)}/{\zeta(t,x)}$.}
\end{equation}
Using (\ref{090504}) and (\ref{090502}), one can check that $q(t,x)$ solves
$$\nu q_t=(\nu q_x)_x +2\nu \frac{\zeta_x}{\zeta}q_x   -c^*q_x,$$
which is a generalization of (\ref{270404bis}).

Recall the function $f(t,x)$ that satisfies the adjoint equation (\ref{0524-18}). A conserved quantity  is
$$I(t)=\int_{c^*t}^\infty \nu(x)f(t,x)p(t,x)dx=\int_{c^*t}^\infty \nu(x)f(t,x)\zeta(t,x)q(t,x)dx.$$
The analog of $I_\alpha$ in (\ref{270410bis}) is
$$I_\alpha(t)=\int_{c^*t}^\infty\eta_\alpha(t,x) \nu(x) p(t,x)dx,$$
with a function $\eta_\alpha(t,x)$ that is exponentially growing as $e^{\alpha(x-c^*t)}$ as $x\to +\infty$. Then
$$\begin{array}{rcl}
\displaystyle\frac{d{I}_\alpha(t)}{dt} & = & \displaystyle\int_{c^*t}^\infty [\pdr{\eta_\alpha}{t}\nu(x)p(t,x)+\eta_\alpha(t,x)((\nu p_x)_x-c^*p_x]dx\\
& = & \displaystyle\int_{c^*t}^\infty \left[\nu(x)\pdr{\eta_\alpha}{t} +\left(\nu(x)\pdr{\eta_\alpha}{x}\right)_x
+c^*\pdr{\eta_\alpha}{x}\right]p(t,x)dx=\mu(\alpha)I_\alpha(t),
\end{array}$$
provided that $\eta_\alpha$ satisfies
$$\nu(x)\pdr{\eta_\alpha}{t} +\left(\nu(x)\pdr{\eta_\alpha}{x}\right)_x
+c^*\pdr{\eta_\alpha}{x}=\mu(\alpha)\nu(x)\eta_\alpha,~~\eta_\alpha(t,c^*t)=0.$$

\begin{lem} \label{lem:etaalphabounds} There is a constant $C > 0$ such that for each $\alpha$ sufficiently small there is a constant $\mu(\alpha)$ and a function $\eta_\alpha(t,x)$ satisfying
\begin{equation}\label{090520}
\nu(x)\pdr{\eta_\alpha}{t} +\left(\nu(x)\pdr{\eta_\alpha}{x}\right)_x
+c^*\pdr{\eta_\alpha}{x}=\mu(\alpha)\nu(x)\eta_\alpha \quad t \in \Rm, \quad x \geq c^* t, \;\;\; t \in \Rm
\end{equation}
and
\[
 \eta_\alpha(t,c^*t)=0, \quad \;\; t  \in \Rm
\]
and 
\[
C \frac{e^{\alpha x} - e^{-\alpha x}}{\alpha} \leq \eta_\alpha(t,x + c^*t) \leq C^{-1} \frac{e^{\alpha x} - e^{-\alpha x}}{\alpha}, \quad \forall \;\;x \geq 0,\;\; t \in \Rm.
\]
In addition, there exists $\mu_0>0$ such that 
\begin{equation}\label{march202}
\hbox{$\mu(\alpha) = \mu_0\alpha^2+O(\alpha^3)$ for all $\alpha > 0$ sufficiently small.}
\end{equation}
\end{lem}

For the homogeneous medium, $\nu(x) \equiv 1$, and the function
\[
\eta_\alpha(t,x) = \frac{e^{\alpha (x - c^*t)} - e^{-\alpha (x-c^*t)}}{\alpha},
\]
satisfies (\ref{090520}) with $\mu(\alpha) = \alpha^2$. 
In the general case, the function $\eta_\alpha$ has exponential asymptotics as $x\to+\infty$
$$\eta_\alpha(t,x)\sim \frac{1}{\alpha} e^{\alpha(x-c^*t)}\bar\eta_\alpha(x),~~\hbox{as $x\to+\infty$,}$$
where $\bar\eta_\alpha(x)$ is a positive periodic solution of
$$\left(\nu(x)\pdr{\bar\eta_\alpha}{x}\right)_x+\alpha(\nu(x)\bar\eta_\alpha)_x+(c^*+\alpha\nu(x))\pdr{\bar\eta_\alpha}{x}
+c^*\alpha(1-\nu(x))\bar\eta_\alpha=(\mu(\alpha)-\alpha^2)\nu(x)\bar\eta_\alpha,$$
and $\mu(\alpha)$ is the corresponding eigenvalue. 

Let us postpone the proof of Lemma \ref{lem:zetabounds} and Lemma \ref{lem:etaalphabounds} and continue with the analysis of $p(t,x)$. We define the second exponential moment by
\[
V_\alpha(t) = \int_{c^*t}^\infty \nu(x) \eta_{2\alpha}(t,x) p(t,x) q(t,x) \,dx =  
\int_{c^*t}^\infty \nu(x) \eta_{2\alpha}(t,x) \zeta(t,x) q^2(t,x) \,dx.
\]
Then
\begin{eqnarray}
&&\frac{d V_\alpha(t)}{dt}=  \int_{c^*t}^\infty \nu (\partial_t  \eta_{2\alpha}) p q \,dx + \int_{c^*t}^\infty \nu \eta_{2\alpha} p_t q \,dx + \int_{c^*t}^\infty \nu \eta_{2\alpha} p q_t \,dx. \no \\
&&~~~~~~~~~
=  \mu(2\alpha) V_\alpha(t) - \int_{c^*t}^\infty \nu (\mathcal{L}^* \eta_{2\alpha}) p q \,dx + \int_{c^*t}^\infty \nu \eta_{2\alpha} p_t q \,dx + \int_{c^*t}^\infty \nu \eta_{2\alpha} p q_t \,dx, \no
\end{eqnarray}
where $\mathcal{L}^* \eta = \nu^{-1}(\nu \eta_x)_x + \nu^{-1} c^* \eta_x$. Since 
\[
p_t = \mathcal{L} p,~~~
q_t = \mathcal{L} q + 2 \frac{\zeta_x}{\zeta} q_x
\]
we have
\be\label{Valph1}\begin{array}{rcl}
V_\alpha'(t) & = & \displaystyle \mu(2\alpha) V_\alpha(t) - \int_{c^*t}^\infty \nu (\mathcal{L}^* \eta_{2\alpha}) p q \,dx 
+ \int_{c^*t}^\infty \nu \eta_{2\alpha} p_t q \,dx + \int_{c^*t}^\infty \nu \eta_{2\alpha} p q_t \,dx\no  \\
& = & \displaystyle  \mu(2\alpha) V_\alpha(t) - \int_{c^*t}^\infty \nu \eta_{2\alpha} \left(p\mathcal{L}q + q \mathcal{L}p\right)  \,dx - 2 \int_{c^*t}^\infty \nu \eta_{2\alpha}  p_x q_x \,dx\\
& & +\displaystyle  \int_{c^*t}^\infty \nu \eta_{2\alpha} p_t q \,dx + \int_{c^*t}^\infty \nu \eta_{2\alpha} p q_t \,dx \no \\
& = &  \displaystyle \mu(2\alpha) V_\alpha(t)  - 2 \int_{c^*t}^\infty \nu \eta_{2\alpha}  p_x q_x \,dx + 2 \int_{c^*t}^\infty \nu \eta_{2\alpha}  p \frac{\zeta_x}{\zeta} q_x \,dx.
\end{array}
\ee
As $p = \zeta q$, we have $p_x = \zeta_x q + \zeta q_x$ and so 
\[
p \frac{\zeta_x}{\zeta} q_x = q \zeta_x q_x = p_x q_x - \zeta (q_x)^2.
\]
Therefore, the last two terms in (\ref{Valph1}) reduce to
\begin{equation}
V_\alpha'(t) = \mu(2\alpha) V_\alpha(t)  - 2 \int_{c^*t}^\infty \nu \eta_{2\alpha}  \zeta (q_x)^2 \,dx. =   \mu(2\alpha) V_\alpha(t) - 2D_\alpha(t) ,\label{Valphp}
\end{equation}
where 
\[
D_\alpha(t) =   \int_{c^*t}^\infty \nu \eta_{2\alpha}  \zeta (q_x)^2 \,dx.
\]

As in the homogeneous case, $V_\alpha(t)$ is the quantity we need to estimate --
we do this by bounding the right side of (\ref{Valphp}).  We claim that there is a constant $C > 0$ such that the inequality
$$D_\alpha(t) \geq C\frac{|V_\alpha(t)|^{5/3}}{|I_\alpha(t)|^{4/3}}$$
holds for all $t > 1$ and $\alpha > 0$ sufficiently small. Since $\nu > 0$ is periodic, this is equivalent to the statement that for any $\alpha > 0$,
\begin{equation}
\left(\int_{c^*t}^\infty \eta_{2\alpha} \zeta q^2 \,dx \right)^{5/3} \leq C \left(\int_{c^*t}^\infty \eta_\alpha \zeta q \,dx\right)^{4/3} \left( \int_{c^*t}^\infty \eta_{2\alpha} \zeta (q_x)^2 \,dx \right). \label{nashsuff}
\end{equation}
By Lemma \ref{lem:zetabounds} and Lemma \ref{lem:etaalphabounds} we may compare the function $\zeta(t,x)$ to the linear function $x - c^*t$, and   $\eta_\alpha(t,x)$ to the function $(e^{\alpha x} - e^{-\alpha x})/\alpha$. That is, for $\alpha > 0$ sufficiently small
\[
 \int_{c^*t}^\infty \eta_{2\alpha} \zeta q^2 \,dx  \leq C_1  \int_{0}^\infty  \frac{e^{2\alpha x} - e^{-2\alpha x}}{2\alpha} x q^2(t,x + c^*t) \,dx  = C_1 \hat V_\alpha,
\]
and 
\[
\int_{c^*t}^\infty \eta_\alpha \zeta q \,dx \geq  C_2 \int_{0}^\infty  \frac{e^{\alpha x} - e^{-\alpha x}}{\alpha} x q(t,x + c^*t) \,dx = C_2 \hat I_\alpha,
\]
and
$$ \int_{c^* t}^\infty \eta_{2\alpha} \zeta (q_x)^2 \,dx \geq C_2 \int_{0}^\infty \frac{e^{2\alpha x} - e^{-2\alpha x}}{2\alpha} x |q_x(t,x+c^*t)|^2 \,dx \geq C_2 \hat D_\alpha.$$
Now (\ref{nashsuff}) follows for all $t > 1$ by applying Lemma \ref{lem:NashR3} with $w(x) = q(t,x + c^* t)$. 

Returning to (\ref{Valphp}) we now have
\[
V_\alpha'(t)  \leq \mu(2\alpha)V_\alpha(t) - C \frac{V_\alpha^{5/3}(t)}{I_\alpha(t)^{4/3}}
\]
where $I_\alpha'(t) = \mu(\alpha)I_\alpha(t)$. For $V_\alpha(t) = e^{\mu(2\alpha)t} Z_\alpha(t)$, this implies the bound
\begin{equation}
Z_\alpha'(t) \leq - C \frac{e^{-t\mu(2\alpha)}e^{ t 5\mu(2\alpha)/3} (Z_\alpha(t))^{5/3} }{e^{t4\mu(\alpha)/3} I_\alpha(0)^{4/3}} = -C\frac{(Z_\alpha(t))^{5/3}}{I_\alpha(0)^{4/3}} e^{t R_\alpha} \label{Zdifalpha}
\end{equation}
for $t \geq 1$, where 
$$R_\alpha  =  \frac{2}{3}\mu(2\alpha) - \frac{4}{3} \mu(\alpha) = \frac{1}{3}\mu(2\alpha)  + O(\alpha^3).$$
We used (\ref{march202}) in the last step above.
We deduce from (\ref{Zdifalpha})  that
\begin{equation}
Z_\alpha(t) \leq C \left( \frac{I_\alpha(0)^{4/3}R_\alpha}{e^{tR_\alpha} - e^{R_\alpha} }\right)^{3/2} =C \frac{I_\alpha(0)^{2}}{(t-1)^{3/2}} \left(  \frac{ t R_\alpha  - R_\alpha}{e^{t R_\alpha} - e^{R_\alpha} }\right)^{3/2}. \label{zRalpha}
\end{equation}
Note that, since $e^{x}$ is a convex function, we have
\[
\frac{b - a}{e^{b} - e^{a}} \leq e^{-a}
\]
for all $b > a$. Moreover, $R_\alpha > 0$ for $\alpha$ sufficiently small, so $R_\alpha t > R_\alpha$ for $t > 1$. 
Hence, (\ref{zRalpha}) implies
\[
Z_\alpha(t) \leq  C \frac{I_\alpha(0)^{2}}{(t-1)^{3/2}} e^{-3 R_\alpha/2} \le  \frac{CI_\alpha(0)^{2}}{(t-1)^{3/2}} .
\]
Therefore, we have
\[
V_\alpha(t) \leq C e^{\mu(2\alpha)t}  \frac{I_\alpha(0)^{2}}{(t-1)^{3/2}},
\]
which is
\[
\left( \int_{c^*t}^\infty \eta_{2\alpha}(t,x) \nu(x) \zeta(t,x) q^2 \,dx \right)^{1/2} \leq 
C \frac{e^{\mu(2\alpha)t}  }{(t-1)^{3/4}}
\int_{0}^\infty \eta_{\alpha}(0,x) \nu(x) \zeta(0,x) p_0(x) \,dx . 
\]
By Lemma \ref{lem:zetabounds} and Lemma \ref{lem:etaalphabounds} and the definition (\ref{march204}) of $ q(t,x)$, this implies
\[
\left( \int_{0}^\infty \frac{e^{2\alpha x} - e^{2\alpha x}}{2\alpha x}  p^2(t,c^*t + x) \,dx \right)^{1/2} 
\leq C \frac{e^{\mu(2\alpha)t}  }{(t-1)^{3/4}} 
\int_{0}^\infty \frac{e^{\alpha x} - e^{\alpha x}}{\alpha} x p_0(x) \,dx.
\]

The rest of the proof of Proposition \ref{prop-0523} now proceeds exactly as in the homogeneous case, in the steps 
following (\ref{080506bis})-(\ref{080508bis}), taking $\alpha=1/\sqrt{t}$, and keeping (\ref{march202}) in mind. 
The only minor technical detail is that the conservation of 
\[
I(t)=\int_0^\infty xp(t,x)dx
\]
is replaced by the conservation of
\[
I(t) = \int_{c^*t}^\infty \nu(x) f(t,x) p(t,x) \,dx,
\]
together with the fact that $m(x - c^*t) \leq f(t,x) \leq m^{-1} (x - c^*t)$ for some $m > 0$, and all $x \geq c^*t$. 
The rest of the argument is essentially identical.~$\Box$

%%%%%%%%%%%%%%%%%%%%%%%%%%%%%%%%%%%%%%%%%

\subsubsection*{Proof of Lemma \ref{lem:zetabounds}}  

The proof is very similar to that of Lemma \ref{lem:fexist}. Recall that the ``linearized" traveling wave is 
\[
\phi(t,x,\lambda) = e^{-\lambda(x - ct)}\psi(x,\lambda).
\] 
At the critical speed, there is another solution of the linearized problem which moves to the right: the function 
\[
\hat \phi(t,x)  = - \frac{\partial}{\partial \lambda} \phi \bigg|_{\lambda = \lambda^*} = 
e^{-\lambda^*(x - c^*t)} \left( (x - c^*t)\psi(x,\lambda^*) - 
\partial_\lambda \psi (x,\lambda^*)\right) = \phi(t,x) \left( (x - c^*t) - \frac{ \psi_\lambda}{\psi}  \right)
\]
is also a solution of the linear equation (\ref{2803-6}). So, if we set 
\[
\Psi(t,x) = (x - c^*t) - \frac{ \psi_\lambda}{\psi} = (x - c^*t) + \chi(x),
\]
(recall (\ref{chidef})) then the two linearized traveling waves are $\phi(t,x)$ and $\hat \phi(t,x) = \phi(t,x) \Psi(t,x)$. 
Therefore, $\Psi$ also satisfies (\ref{090504}):
\[
\nu (x)\pdr{\Psi}{t} = \left(\nu(x)\pdr{\Psi}{x}\right)_x-c^*\pdr{\Psi}{x}.
\]
In the homogeneous case we have $\psi(x,\lambda)\equiv 1$, hence we take $\zeta(t,x)=\Psi(t,x) = (x - c^*t)$. In general, however, $\Psi$ doesn't satisfy the Dirichlet boundary condition at $x=c^*t$ in (\ref{zetabc}), therefore we can not take $\zeta(t,x)=\Psi(t,x)$. Instead, we take $\zeta(t,x)$ to be the limit (as $n \to \infty$) of a sequence of functions $\{\zeta^{(n)}(t,x) \}_{n=1}^\infty$ which satisfy
\begin{eqnarray}
&& \nu(x) \zeta^{(n)}_t =(\nu(x) \zeta^{(n)}_x)_x - c^* \zeta^{(n)}_x, \quad x > c^* t, \;\; t \geq -n \no \\
&&\zeta^{(n)}(t,c^*t) = 0, \quad t \geq -n \no \\
&&\zeta^{(n)}(-n,x) = \max(0,\Psi(-n,x)), \quad x \geq -c^*n.\no
\end{eqnarray}
The maximum principle implies that for some constant $C$,
$$\Psi(t,x) - C \leq \zeta^{(n)}(t,x) \leq \Psi(t,x) + C, $$
holds for all $t \geq -n$ and $x \geq c^* t$. The rest follows as in the proof of Lemma \ref{lem:fexist}. \qed

%%%%%%%%%%%%%%%%%%%%%%%%%%%%%%%%%%%%%%%%%

\subsubsection*{Proof of Lemma \ref{lem:etaalphabounds}}

{\bf The eigenvalue asymptotics for $\alpha\ll 1$.} 
The only remaining  ingredient in the proof of Proposition~\ref{prop-0523} is the
proof of Lemma \ref{lem:etaalphabounds}. First, we prove the asymptotics (\ref{march202}) for $\mu(\alpha)$.
Consider the periodic eigenvalue problem 
\be\label{090526}\left\{\begin{array}{l}
\displaystyle\left(\nu(x)\pdr{\eta}{x}\right)_x+\alpha(\nu(x)\eta)_x+(c^*+\alpha\nu(x))\pdr{\eta}{x}
+c^*\alpha(1-\nu(x))\eta=\gamma(\alpha) \nu(x)\eta,\no\\
\eta(x + 1) = \eta(x) > 0,
\end{array}\right.
\ee
with $\gamma(\alpha) = \mu(\alpha) - \alpha^2$ and the normalization 
$$\int_0^1 \nu(x)  \eta(x) \,dx = 1.$$
Observe that $\gamma(0) = 0$ and $\eta(x, \alpha = 0) \equiv 1$. Moreover,  as $\gamma(0) = 0$ is a
simple eigenvalue, $\gamma(\alpha)$ is an
 analytic function of $\alpha$, for $\alpha$ sufficiently small.  The function $\eta' = {\partial\eta}/{\partial\alpha}$ satisfies
\[
\left(\nu(x)\pdr{\eta'}{x}\right)_x+\alpha(\nu(x)\eta')_x+(c^*+\alpha\nu(x))\pdr{\eta'}{x}
+c^*\alpha(1-\nu(x))\eta' + (\nu\eta)_x + \nu \pdr{\eta}{x} + c^* (1 - \nu) \eta=\gamma\nu\eta' + \gamma' \nu \eta.  
\]
Setting $\alpha = 0$ we obtain:
\begin{eqnarray}
&& \!\!\!\!\!\left(\nu(x)\pdr{\eta'}{x}\right)_x + c^*\pdr{\eta'}{x}
+ \nu_x + c^* (1 - \nu) = \gamma' \nu . \label{etaprimezer}
\end{eqnarray}
Integrating (\ref{etaprimezer}), we conclude that $\gamma'(0) = 0$. Next, $\eta''$ solves
\begin{eqnarray*}
&& \!\!\!\!\!\!\!\!\!\!
\left(\nu(x)\pdr{\eta''}{x}\right)_x+\alpha(\nu(x)\eta'')_x+(c^*+\alpha\nu(x))\pdr{\eta''}{x}
+c^*\alpha(1-\nu(x))\eta'' + 2(\nu\eta')_x + 2\nu \pdr{\eta'}{x} + 2c^* (1 - \nu) \eta' \no \\
& & \quad \quad \quad \quad  \quad \quad \quad   \quad \quad \quad \quad \quad \quad \quad =\gamma \nu\eta'' + 2\gamma' \nu\eta' + \gamma'' \nu \eta.
\end{eqnarray*}
So, at $\alpha = 0$ we have
\begin{eqnarray}
 \left(\nu(x)\pdr{\eta''}{x}\right)_x+c^*\pdr{\eta''}{x} + 2(\nu\eta')_x + 2\nu \pdr{\eta'}{x} + 2c^* (1 - \nu) \eta'  =  \gamma'' \nu. \no
\end{eqnarray}
Integrating this equation, we obtain
\begin{equation}
\gamma'' = 2 \int_0^1 \left( \nu \pdr{\eta'}{x} + c^*(1 - \nu)\eta' \right) \,dx. \label{gppeqn}
\end{equation}
Since $\gamma'(0) = 0$, (\ref{etaprimezer}) implies that
\[
c^*(1 - \nu) = - \nu_x - c^* \pdr{\eta'}{x} - \left( \nu \pdr{\eta'}{x} \right)_x.
\]
Plugging this into (\ref{gppeqn}), we obtain
\[
\gamma''(0) = 4\int_0^1  \nu \pdr{\eta'}{x} \,dx + 2 \int_0^1 \nu \left(\pdr{\eta'}{x}\right)^2 \,dx.
\]
Since $4 y + 2y^2 \geq -2$ for all $y \in \Rm$, we conclude that
\[
\gamma''(0) \geq - 2 \int_0^1 \nu(x) \,dx = -2
\]
with equality if and only if $\pdr{\eta'}{x} \equiv -1$.  Since $\eta'$ is periodic, $\pdr{\eta'}{x} = -1$ cannot hold at all $x$, so we must have $\gamma''(0) > -2$. Finally, since $\mu(\alpha) = \alpha^2 + \gamma(\alpha)$, we have $\mu''(0) = 2 + \gamma''(0) > 0$,
proving (\ref{march202}).

Let us now denote the eigenfunction of (\ref{090526}) by $\bar\eta_\alpha$ to indicate its dependence on $\alpha$.

\begin{cor} \label{cor:betaalpha}
There is a constant $C$ such that for all $\alpha > 0$ sufficiently small, there is
$\beta(\alpha) > 0$ with $\mu(-\beta) = \mu(\alpha)$ and such that
\begin{equation}
\left| \frac{\beta}{\alpha} - 1 \right| \leq C \alpha, \label{betaalphasymp}
\end{equation}
and
\[
\sup_{x} \left|  \bar \eta_\alpha(x) - 1\right| \leq C \alpha,
~~~~\sup_{x} \left|  \bar \eta_\beta(x) - 1\right| \leq C \alpha.
\]
\end{cor}

\noindent{\bf Proof.} The existence of such a $\beta$ satisfying (\ref{betaalphasymp}) 
follows from the fact that $\mu(\alpha) \sim C \alpha^2$ for $\alpha$ small. The bounds on $\bar \eta_\alpha$ and 
$\bar \eta_\beta$ follow from elliptic regularity and the fact that for $\alpha = 0$, $\bar \eta_{0}(x) \equiv 1$. \qed\hfill\break

\noindent{\bf Construction of the function $\eta_\alpha(t,x)$.}
Continuing with the proof of Lemma \ref{lem:etaalphabounds}, choose $\beta = \beta(\alpha) > 0$ 
according to Corollary \ref{cor:betaalpha} and consider the terminal value problem 
\begin{equation} \label{etaterm}
\nu(x)\pdr{\eta_{\alpha,T}}{t} +\left(\nu(x)\pdr{\eta_{\alpha,T}}{x}\right)_x
+c^*\pdr{\eta_{\alpha,T}}{x}=\mu(\alpha)\nu(x)\eta_{\alpha,T} \quad t < T, \quad x \geq c^* t
\end{equation}
with the terminal condition $\eta_{\alpha,T}(T,x) \geq 0$ to be determined. The function $\eta_\alpha(t,x)$ of 
Lemma \ref{lem:etaalphabounds} will be defined as $\lim_{T \to \infty} \eta_{\alpha, T}(t,x)$.  
Observe that for any constant $C$, the function
\[
\alpha^{-1} e^{\alpha (x - c^*t)} \bar \eta_{\alpha}(x) - C \beta^{-1}e^{- \beta(x - c^*t)}\bar \eta_{\beta}(x)
\]
satisfies (\ref{etaterm}), since $\mu(-\beta) = \mu(\alpha)$. If we choose the constant  
\[
C_u = \frac{\beta}{\alpha} \min_{x} \frac{\bar \eta_\alpha(x)}{\bar \eta_{\beta}(x)} > 0,
\]
then the function
\[
h_u(t,x) =  \alpha^{-1} e^{\alpha (x - c^*t)} \bar \eta_{\alpha}(x) - C_u \beta^{-1}e^{- \beta(x - c^*t)}\bar \eta_{\beta}(x)
\]
satisfies $h_u(t,c^*t) \geq 0$ for all $t \in \Rm$. Similarly, if we choose
\[
C_l = \frac{\beta}{\alpha} \max_{x} \frac{\bar \eta_\alpha(x)}{\bar \eta_{\beta}(x)} > 0,
\]
then the function
\be
h_l(t,x) =  \alpha^{-1} e^{\alpha (x - c^*t)} \bar \eta_{\alpha}(x) - C_l \beta^{-1} e^{- \beta(x - c^*t)}\bar \eta_{\beta}(x)  \label{hldef}
\ee
satisfies $h_l(t,c^*t) \leq 0$ for all $t \in \Rm$. Now, if we choose the terminal condition for $\eta_{\alpha,T}$ to be
\[
\eta_{\alpha,T}(T,x) = \max \left( 0 ,h_\ell(T,x) \right), 
\]
the maximum principle implies that 
\begin{equation}
h_l(t,x) \leq \eta_\alpha(t,x) \leq h_u(t,x) \label{huhl}
\end{equation}
holds for all $t \leq T$ and $x \geq c^*t$. Although the constants $C_u$ and $C_l$ depend on $\alpha$, 
Corollary \ref{cor:betaalpha} implies that
\[
C_u = 1 + O(\alpha) \quad \quad \text{and}\quad \quad C_l = 1 + O(\alpha)
\]
as $\alpha \to 0$.

Now, we claim there are constants $L > 0$ and $M > 0$, independent of $T$, such that
\begin{equation}
M \frac{e^{\alpha x} - e^{-\alpha x}}{\alpha}\leq \eta_{\alpha,T}(t,x + c^*t) \leq M^{-1} \frac{e^{\alpha x} - e^{-\alpha x}}{\alpha}, 
\label{etaclaim}
\end{equation}
for all $x > L$ and $t \leq T$, and all $\alpha$ sufficiently small. 
Given this claim, parabolic regularity and the maximum principle imply that there is a constant $b > 0$ (also independent of $T$) such that
\[
b < \frac{\partial \eta_{\alpha,T}}{\partial x}(t,c^*t) < b^{-1}
\]
holds for all $t \leq T -1$ and $\alpha > 0$ sufficiently small. Since 
\[
\frac{d}{dx} \left(\frac{e^{\alpha x} - e^{-\alpha x}}{\alpha } \right)\bigg|_{x = 0} = 2,
\]
it follows, by parabolic regularity, that
$$C \frac{e^{\alpha x} - e^{-\alpha x}}{\alpha}\leq \eta_{\alpha,T}(t,x + c^*t) \leq C^{-1} \frac{e^{\alpha x} - e^{-\alpha x}}{\alpha},$$
for all $x \geq 0$ and $t \leq T -1$, with a constant $C$ that is independent of $T$. 
Then letting $T \to +\infty$ we may take a subsequence of functions $\eta_{\alpha,T_k}(x,t) $ such that 
$T_k \to \infty$ and $\eta_{\alpha,T_k}$ converges locally uniformly to a function $\eta_{\alpha}(t,x)$ 
satisfying all the criteria of Lemma \ref{lem:etaalphabounds}.\hfill\break

\noindent{\bf The proof of (\ref{etaclaim}).} Let us derive the upper bound in (\ref{etaclaim}). Because of (\ref{huhl}), it suffices to show that 
\begin{equation}
h_u(t,x + c^*t) \leq M^{-1} \frac{e^{\alpha x} - e^{-\alpha x}}{\alpha} \label{husuff}
\end{equation}
holds for all $t \in \Rm$ and $x \geq L$, with $L > 0$ and $M$ being independent of $\alpha$. Let us write $h_u(t,x)$ as
\[
h_u(t,x + c^*t) =  \alpha^{-1}\bar \eta_{\alpha}(x + c^*t) \left( e^{\alpha x}  - 
C_u \frac{\alpha}{\beta} \frac{\bar \eta_{\beta}(x + c^*t)}{\bar \eta_\alpha(x + c^*t)}e^{- \beta x} \right). 
\]
Therefore, since $\bar \eta_{\alpha}$ is uniformly bounded in $x$, independently of $\alpha \in (0,1)$, the upper bound (\ref{husuff}) holds if
\[
\left( e^{\alpha x}  - C_u \frac{\alpha}{\beta} \frac{\bar \eta_{\beta}(x + c^*t)}{\bar \eta_\alpha(x + c^*t)}e^{- \beta 	x} \right) \leq M_{2}\left(e^{\alpha x} - e^{-\alpha x}\right)
\]
for some constant $M_2$, which is equivalent to
\begin{equation}
e^{-2\alpha x}  \left( M_2  - C_u \frac{\alpha}{\beta} \frac{\bar \eta_{\beta}(x + c^*t)}{\bar \eta_\alpha(x + c^*t)} e^{- (\beta - \alpha) x} \right) \leq M_{2} - 1. \label{alphequiv}
\end{equation}
Since $C_u$, $\bar \eta_\alpha$, $\bar \eta_{\beta}$ are positive, this inequality certainly holds if 
\[
e^{-2 \alpha x} M_2 \leq M_2 - 1.
\]
So, if we set $M_2 = 2$, then (\ref{alphequiv}) holds for all $x \geq  \ln(2)/(2\alpha)$. Now consider (\ref{alphequiv}) for $x \leq  \ln(2)/(2\alpha)$. By Corollary \ref{cor:betaalpha}
\[
C_u \frac{\alpha}{\beta} \frac{\bar \eta_{\beta}(x + c^*t)}{\bar \eta_\alpha(x + c^*t)} = 1 + O(\alpha)
\]
as $\alpha \to 0$, uniformly in $x$ and $t$. Moreover, $\beta - \alpha = O(\alpha^2)$, so that for $x \leq  \ln(2)/(2\alpha)$, we have 
\[
C_u \frac{\alpha}{\beta} \frac{\bar \eta_{\beta}(x + c^*t)}{\bar \eta_\alpha(x + c^*t)}  e^{- (\beta - \alpha) x} = (1 + O(\alpha)) 
\]
Therefore, with $M_2 = 2$ and $x \leq  \ln(2)/(2\alpha)$, inequality (\ref{alphequiv}) becomes
$$e^{-2\alpha x} \leq \frac{M_2 - 1}{ M_2 - 1 + O(\alpha)} = 1 - O(\alpha).$$
Hence there is a constant $L$ such that (\ref{alphequiv}) holds for all $x \geq L$ and $t \in \Rm$, and all $\alpha$ sufficiently small. This establishes the upper bounds in (\ref{husuff}) and (\ref{etaclaim}).

In a similar manner, we now we prove the lower bound in (\ref{etaclaim}).  It suffices to show that 
\begin{equation}
h_l(t,x) \geq M \frac{e^{\alpha x} - e^{-\alpha x}}{\alpha} \label{hsufflower}
\end{equation}
holds for all $t \in \Rm$ and $x \geq L$.  Let us write $h_l(t,x)$ as
\[
h_l(t,x + c^*t) =  \alpha^{-1}\bar \eta_{\alpha}(x + c^*t) \left( e^{\alpha x}  - C_l \frac{\alpha}{\beta} \frac{\bar \eta_{\beta}(x + c^*t)}{\bar \eta_\alpha(x + c^*t)}e^{- \beta x} \right). 
\]
Therefore, since $\bar \eta_{\alpha}(x)$ is uniformly bounded away from zero, independently of $\alpha \in (0,1)$, the lower bound (\ref{hsufflower}) holds if
\[
M_3 \left( e^{\alpha x}  - C_l \frac{\alpha}{\beta} \frac{\bar \eta_{\beta}(x + c^*t)}{\bar \eta_\alpha(x + c^*t)}e^{- \beta 	x} \right) \geq e^{\alpha x} - e^{-\alpha x}
\]
for some constant $M_3$, which is equivalent to
\begin{equation}
M_3 - 1 \geq  M_3 C_l \frac{\alpha}{\beta} \frac{\bar \eta_{\beta}(x + c^*t)}{\bar \eta_\alpha(x + c^*t)} e^{- (\beta + \alpha) x} - e^{-2\alpha x}. \label{alphequiv2}
\end{equation}
This bound certainly holds if 
\[
M_3 - 1 \geq M_3 C_l \frac{\alpha}{\beta} \frac{\bar \eta_{\beta}(x + c^*t)}{\bar \eta_\alpha(x + c^*t)} e^{- (\beta + \alpha) x}.
\]
By Corollary \ref{cor:betaalpha} we know that 
\[
C_l \frac{\alpha}{\beta} \frac{\bar \eta_{\beta}(x + c^*t)}{\bar \eta_\alpha(x + c^*t)} = 1 + O(\alpha) \leq 2
\]
uniformly in $x$ and $t$, if $\alpha$ is sufficiently small. So, if we set $M_3 = 2$, then (\ref{alphequiv2}) holds for all $x \geq  \ln(2)/\alpha$.

Now consider (\ref{alphequiv2}) for $x \leq  \ln(2)/\alpha$. Recall that, $\beta + \alpha = 2\alpha + O(\alpha^2)$, so that for $x \leq  \ln(2)/\alpha$, we have 
\[
C_l \frac{\alpha}{\beta} \frac{\bar \eta_{\beta}(x + c^*t)}{\bar \eta_\alpha(x + c^*t)}  e^{- (\beta - \alpha) x} = e^{-2\alpha x}(1 + O(\alpha)). 
\]
Therefore, with $M_3 = 2$ and $x \leq  \ln(2)/\alpha$, inequality (\ref{alphequiv2}) becomes
$$M_3 - 1 \geq \left( M_3(1 + O(\alpha)) - 1 \right) e^{-2\alpha x} ,$$
which is
\[
e^{-2\alpha x} \leq \frac{1}{2(1 + O(\alpha)) - 1} = 1 - O(\alpha).
\]
Hence there is a constant $L$ such that (\ref{alphequiv2}) holds for all $x \geq L$ and $t \in \Rm$, and all $\alpha$ sufficiently small. This proves the lower bound in (\ref{hsufflower}) and in (\ref{etaclaim}). This completes the proof of Lemma~\ref{lem:etaalphabounds},
and the proof of Propositions~\ref{cor-0524} and \ref{prop-0523} is also now complete.~\qed

%%%%%%%%%%%%%%%%%%%%%%%%%%%%%%%%%%%%%%%%%
%%%%%%%%%%%%%%%%%%%%%%%%%%%%%%%%%%%%%%%%%

\section{Proof of Proposition~\ref{thm-plan1}}

Recall that the upper bound in~Theorem~\ref{thm-delay-apr5} was reduced in Section~\ref{sec:upper} to the
proof of Proposition~~\ref{thm-plan1} that we present in this section.
Let $\tilde p(\tau,x)$ be as in this proposition, that is
\begin{equation}\label{march208}
\left(1-\omega(\tau)\right)\tilde p_\tau = \tilde p_{xx} + 2 \frac{\phi_x}{\phi} \tilde p_x,~~x\ge c^*\tau,
\end{equation}
with the Dirichlet boundary condition $\tilde p(\tau,c^*\tau)=0$. The coefficient $\omega(\tau)$ satisfies  
$\omega(\tau) \sim {3}/(2{c^* \tau})$ as $\tau \to \infty$, and  $|\omega(\tau)|\le{C}/{\tau}$,  $|\omega'(\tau)|\le C/\tau^2$ for $\tau>\tau_0$.  
The general philosophy is that the correction $\omega(\tau)$ does not play a role in most of the decay estimates, and the function
$\tilde p(t,x)$ behaves essentially as $p(t,x)$, which is the solution of (\ref{march208}) with $\omega(\tau)=0$, and which we have
studied in detail in the preceding sections.
We will need the following steps to prove Proposition~\ref{thm-plan1}. The key step is to establish that $\tilde p(t,x)$ decays as $C/\tau$ at 
positions of the order $c^*\tau+O(\sqrt{\tau})$.

\begin{prop}\label{plan-prop1}
For any $L_0>0$ and $\varepsilon > 0$, there exists $C_\eps>0$ so that 
$$\farc{1}{C_\eps \tau}\le \tilde p(\tau,c^*\tau+L_0+\eps\sqrt{\tau})\le\frac{C_\eps}{\tau}$$
holds for all $\tau \geq 1$.
\end{prop}

This is a direct generalization of  Proposition \ref{cor-0524} and Lemma \ref{prop-0504-2bis}  to the case $\omega(\tau)\neq 0$.
We will also need a more or less explicit solution of the approximate equation that we will need compare to $\tilde p(t,x)$.
It is described in the next proposition.

\begin{prop}\label{plan-prop2}
Let $\bar \chi \in \Rm$ and let $\chi(x)$ be as in (\ref{chidef}). There is a function $\theta^{app}(\tau,x)$ such that for any 
$\sigma > 0$, $\theta^{app}(\tau,x)$ satisfies
$$(1-\omega(\tau))\pdr{\theta^{app}}{\tau} - \theta_{xx}^{app} - 2 \frac{\phi_x}{\phi} \theta_x^{app} =O(\tau^{-3}),~~c^*\tau<x<c^*\tau+ \sigma \sqrt\tau, \;\; \tau \geq 1,$$
and there is a constant $C$ (depending on $\sigma$ and $m$) such that
\be\label{thetaappbound}
\left| \theta^{app}(\tau,x) - \frac{x - c^* \tau + \chi(x) + \bar \chi}{\tau^{3/2}}\; e^{-\frac{(x - c^* \tau)^2}{4(1 + \kappa)\tau}} \right| \leq C \tau^{-3/2} \left( \frac{x - c^* \tau}{\sqrt{\tau}} \right)^2 + O(\tau^{-2})
\ee
holds for all $x \in [c^* \tau, c^* \tau + \sigma \sqrt{\tau}]$ and $\tau \geq 1$. The constant $\kappa$ in the exponential factor is defined by formula~(\ref{kappadef}) below and satisfies $1 + \kappa > 0$. 
\end{prop}

The following refinement of the approximate solution satisfies the exact problem.

\begin{prop}\label{plan-prop3}
Let $\sigma > 0$ be fixed, and let $\theta^{app}(\tau,x)$ be defined as in Proposition \ref{plan-prop2} for some $\bar \chi \in \Rm$. Let $\xi(\tau,x)$ solve 
\begin{equation}\label{plan-eq28}
(1-\omega(\tau))\pdr{\xi}{\tau}= \xi_{xx} + 2 \frac{\phi_x}{\phi} \xi_x,~~x \in (c^*\tau , c^* \tau + \sigma \sqrt{\tau}), \;\; \tau > 1
\end{equation}
with the boundary conditions
\begin{eqnarray}\label{plan-eq30}
\xi(\tau,c^*\tau) & = & \theta^{app}(\tau,c^* \tau ), \no \\
\xi(\tau,c^*\tau+ \sigma \sqrt{\tau}) & = & \theta^{app}(\tau,c^*\tau+ \sigma \sqrt{\tau}).
\end{eqnarray}
Then there is $\tau_0 > 0$ such that
$$|\xi(\tau,x)-\theta^{app}(\tau,x)|\le\frac{C}{\tau^{3/2}},~~c^*\tau < x <c^*\tau+ \sigma \sqrt{\tau}$$
holds for all $\tau \geq \tau_0$. 
\end{prop}

Observe that by choosing $\bar \chi > \norm{\chi}_\infty$ in Proposition \ref{plan-prop2}, we may arrange that 
$\theta^{app}(\tau,c^*\tau) > 0$ for $\tau$ suffiently large. Similarly, with $\bar \chi < -\norm{\chi}_\infty$, we
 have $\theta^{app}(\tau,c^*\tau) < 0$ for $\tau$ sufficiently large. Let us define $\theta^{app}_+$ to be a solution 
 with $\bar \chi = 2 \norm{\chi}_\infty$ and $\theta^{app}_+(\tau,c^*\tau) > 0$; let $\theta^{app}_-$ to be a solution 
 with $m = -2 \norm{\chi}_\infty$ and $\theta^{app}_-(\tau,c^*\tau) < 0$. To prove Proposition~\ref{thm-plan1}, 
 we wish to compare $\tilde p(\tau,x)$ with the functions $\theta^{app}_\pm$. We know from Proposition \ref{plan-prop2} that
\[
\left|\theta^{app}_\pm(\tau,c^*\tau+ \sigma \sqrt{\tau})-\farc{C\sigma}{\tau}\right|\le\farc{C}{\tau^{3/2}}.
\]
Combining this with Proposition~\ref{plan-prop1}, we see that there must be $C_1 > 0$ such that
\[
 \tilde p(\tau, c^*\tau+\sigma  \sqrt{\tau}) \leq C_1 \theta^{app}_+(\tau,c^*\tau+\sigma  \sqrt{\tau})
\]
\[
 \tilde p(\tau, c^*\tau+\sigma  \sqrt{\tau}) \geq C_1^{-1} \theta^{app}_-(\tau,c^*\tau+\sigma  \sqrt{\tau}),
\]
for all $\tau \geq 1$. Now if $\xi_\pm(\tau,x)$ solve (\ref{plan-eq28}) for $\tau \geq 1$ with 
the boundary conditions (\ref{plan-eq30}) using $\theta^{app} = \theta^{app}_\pm$, we have
\begin{eqnarray}
&& \xi_+(\tau,c^* \tau) = \theta^{app}_+ (\tau,c^* \tau) > 0 = C_1^{-1} \tilde p (\tau,c^* \tau) \no \\
&& \xi_+(\tau,c^*\tau+ \sigma\sqrt{\tau}) = \theta^{app}_+ (\tau,c^*\tau+ \sigma\sqrt{\tau}) \geq C_1^{-1} \tilde p (\tau,c^* \tau)\no
\end{eqnarray}
and
\begin{eqnarray}
\xi_-(\tau,c^* \tau) = \theta^{app}_- (\tau,c^* \tau) < 0 = C_1 \tilde p (\tau,c^* \tau) \no \\
\xi_-(\tau,c^*\tau+ \sigma\sqrt{\tau}) = \theta^{app}_- (\tau,c^*\tau+ \sigma \sqrt{\tau}) \leq C_1 \tilde p (\tau,c^* \tau)\no
\end{eqnarray}
The maximum principle implies
$$C_1^{-1} \xi_-(\tau,x)\le \tilde p(\tau,x)\le C_1 \xi^+(\tau,x),$$
holds for all $\tau $ sufficiently large and $x \in [c^* \tau, c^* \tau +  \sigma \sqrt{\tau}]$.

Proposition~\ref{plan-prop3} implies that for any $\delta>0$ there exists
$x_\delta$ so that
$$|\xi_\pm(\tau,x)-\theta^{app}_\pm(\tau,x)|\le\delta \theta^{app}_\pm(\tau,x),~~c^*\tau+x_\delta<x<c^*\tau+\eps{\sqrt{\tau}},$$
if $\tau \geq \tau_0$. It follows that
$$\frac{C_1^{-1}}{2}\theta^{app}_-(\tau,x)\le \tilde p(\tau,x)\le 2{C_1}\theta^{app}_+(\tau,x),
~~c^*\tau+x_\delta<x<c^*\tau+\eps{\sqrt{\tau}},$$
for all $\tau \geq \tau_0$. In view of (\ref{thetaappbound}) and parabolic regularity, 
the conclusion of Proposition~\ref{thm-plan1} now follows.~\qed
 
%%%%%%%%%%%%%%%%%%%%%%%%%%%%%%%%%%%%%%%%%
 
\subsection*{The proof of Proposition~\ref{plan-prop1}}

The proof of Proposition~\ref{plan-prop1} is as in the case $\omega(\tau)=0$ 
(i.e. Proposition \ref{cor-0524} and Lemma \ref{prop-0504-2bis}) but a little more technical -- we focus only on the differences. The first ingredient needed is a quantity that is bounded from above and below.

\begin{lem}
\label{l5.1}
Let $\tilde p(\tau,x)$ be as in Proposition \ref{thm-plan1}. There is $C>0$ such that 
$$
C^{-1} \leq \int_{c^* \tau}^{+\infty}(x - c^* \tau) \tilde p(\tau,x)\,dx \leq C, \quad \forall \; \tau \geq 0.
$$
\end{lem}

\noindent{\bf Proof.} It suffices to bound the integral
$$
I(\tau)=\int_{c^*\tau}^{+\infty}\nu (x) (1-\omega(\tau))f(\tau,x)\tilde p(\tau,x)\,dx,
$$
where $f(\tau,x)$ is the function defined in Lemma \ref{lem:fexist}, with $m (x - c^* \tau) \leq f \leq m^{-1} (x - c^* \tau)$. 
In the case $\omega \equiv 0$, $I(\tau)$ is conserved. We compute:
\begin{eqnarray}
\frac{dI}{d\tau}  = -\omega'\int_{c^*\tau}^{+\infty}\nu f\tilde p \, dx - \omega\int_{c^*\tau}^{+\infty} \nu f_\tau \tilde p\,dx 
 =  O(\tau^{-2}) I(\tau) - \omega\int_{c^*\tau}^{+\infty} \nu f_\tau \tilde p\,dx. \label{IprimeO}
\end{eqnarray}
For an upper bound on $I(\tau)$, we treat the spurious term $\di\int_{c^*\tau}^{+\infty}\nu f_\tau \tilde p \,dx$ as follows:
$$
\di\int_{c^*\tau}^{+\infty}\nu f_\tau \tilde p \, dx=\int_{c^*\tau}^{c^*\tau+\tau^{1/4}}\nu f_\tau \tilde p \,dx +\int_{c^*\tau+\tau^{1/4}}^{+\infty} \nu f_\tau \tilde p \,dx :=I\!I+I\!I\!I.
$$
By parabolic regularity, there is a constant $C>0$ such that $\vert\partial_\tau f(\tau,x)\vert\leq C$, hence
$$
|I\!I\!I|\leq C \tau^{-1/4}\int_{c^*\tau}^{+\infty}x \tilde p \,dx\leq C\tau^{-1/4}\int_{c^*\tau}^{+\infty}\nu(x) f(\tau,x) \tilde p(\tau,x)\,dx.
$$

Recall that equation (\ref{plan-eq2}) for $\tilde p$ is equivalent to
\be
(1 - \omega(\tau)) \tilde p_\tau = \frac{1}{\nu(x)} ( \nu (x) \tilde p_x)_x - \frac{c^*}{\nu(x)} \tilde p_x, \quad x > c^* \tau \label{qselfadj}
\ee
with $\tilde p(\tau,c^* \tau) =0$. A simple time change  so that
\[
d\tau'=\frac{d\tau}{1-\omega(\tau)},
\]
shows the heat kernel bounds of ~\cite{Norris} in the whole space
hold (with the time change) for the perturbed equation
\[
\left(1-\omega(\tau)\right)P_\tau =  \frac{1}{\nu(x)} ( \nu (x) P_x)_x - \frac{c^*}{\nu(x)} P_x ,~~x\in\mathbb{R}.
\]
In particular, we have
$$
|P(\tau,x)|\leq C\tau^{-1/2}\int_{\mathbb{R}} |P(0,y)|dy. 
$$
So, because $\tilde p(\tau,x)$ is less than the solution of (\ref{qselfadj}) in the whole space with the same 
initial data $\tilde p(0,\cdot)$, we have:
$$
|I\!I|\leq  C \tau^{-1/2} \int_{c^* \tau}^{c^* \tau + \tau^{1/4}}\int_{\mathbb{R}} |\tilde p(0,y)|\, dy \, dx =
C\tau^{-1/4} \int_{0}^{+\infty}\tilde p(0,x) \,dx.
$$
Gathering these estimates we conclude
$$
I'(\tau) \leq   O(\tau^{-2})I + O(\tau^{-5/4})I + O(\tau^{-5/4}),
$$
which implies the existence of $C>0$ such that $I(\tau)\leq C(1 + I(0))$.

For a lower bound, note that $f_\tau \leq 0$, while $\nu, \tilde p \geq 0$. Therefore, the term 
\[
- \omega\di\int_{c^*\tau}^{+\infty}\nu f_\tau \tilde p \,dx
\]
in (\ref{IprimeO}) is non-negative. This implies $I'(\tau) \geq O(\tau^{-2}) I$, so that $I(t) \geq C I(0) > 0$,
with some constant $C > 0$.  ~$\Box$\hfill\break

The main step in the proof of Proposition~\ref{plan-prop1} is an estimate on the quantity
$$
V_\alpha(\tau)=(1-\omega(\tau))\int_{c^*\tau}^{+\infty}\nu(x) \eta_{2\alpha}(\tau,x)\tilde p(\tau,x)q(\tau,x)dx,
$$
where 
$$q(\tau,x) = \frac{\tilde p(\tau,x)}{\zeta(\tau,x)}$$
and $\zeta(\tau,x)$ is defined by Lemma \ref{lem:zetabounds}; it solves the unperturbed equation
(i.e. with $\omega=0$) and grows linearly. The function $\eta_\alpha(\tau,x)$ is defined by Lemma \ref{lem:etaalphabounds}. 
The time derivatives of the   functions $\eta$ and $\zeta$ will have to be examined, and this is the object of the following

\begin{lem} \label{l6.1} 
(i). There is a constant $C>0$ such that $\vert\partial_\tau\zeta(\tau,x)\vert\leq C$ for all $x > c^* \tau$.\hfill\break
\noindent (ii). There is a consant $C$ such that $\vert\partial_\tau\eta_\alpha(\tau,x)\vert\leq C$ for all $x \in (c^* \tau, c^*\tau + \alpha^{-1})$.\hfill\break
\noindent (iii). There is a constant $C$ such that $\vert\partial_\tau\eta_\alpha(\tau,x)\vert\leq C\alpha\eta_\alpha(\tau,x)$ for all $x > c^* \tau$.
\end{lem}

\noindent{\bf Proof.} 
Part (i) just comes from parabolic regularity. As for Part (ii), we come back to the notations of Lemma \ref{lem:etaalphabounds}.
Consider $T>0$, at $\tau=T$ we have, just using the equation for $\eta_\alpha$:
$$
\partial_\tau \eta_{\alpha}(T,x)=O(e^{\alpha(x-c^*T)}+e^{-\alpha(x-c^*T)})+d\mu_\alpha(x)
$$
where $\mu_\alpha$ is a measure carried by the (compact) zero set of the function $h_l$, which was defined at (\ref{hldef}), and whose mass is uniformly bounded with respect to $\alpha$.  So, applying the equation for $\partial_\tau \eta_{\alpha}$ - recall that it solves the same equation as $\eta_\alpha$:
$$
\partial_\tau \eta_{\alpha}(T-1,x)=O(e^{\alpha(x-c^*T)}+e^{-\alpha(x-c^*T)})+O(1)=O(e^{\alpha(x-c^*T)}).
$$
Running the equation for $\tau\leq T-1$ yields 
$$
\vert\partial_\tau \eta_{\alpha}(\tau,x)\vert\leq Ce^{\alpha(x-c^*\tau)}\bar\eta_\alpha(x),
$$
and so $\partial_\tau \eta_{\alpha}(\tau,x)=O(e^{\alpha(x-c^*\tau)})$, which is sufficient to prove the claim. \qed

\vspace{0.2in}

\noindent{\bf Proof of Proposition \ref{plan-prop1}.} A straightforward computation shows that 
\begin{eqnarray}
\frac{dV_\alpha}{d\tau} & = & (\mu(2 \alpha) - \omega') V_\alpha  - \omega \int \nu \eta_\tau \tilde p q \,dx - 2\int \nu \eta \tilde p_x q_x \,dx \no \\
& & + 2 \int \nu \eta \tilde p  \frac{\zeta_x}{\zeta} q_x \,dx + \omega \int \nu \eta \tilde p q  \frac{\zeta_\tau}{\zeta}  \,dx \no \\
& = &
 (\mu(2 \alpha) - \omega') V_\alpha  - 2D_\alpha  + \omega \int \nu \left( \eta \tilde p q \frac{\zeta_\tau}{\zeta} - \eta_\tau \tilde p q \right) \,dx \no
\\
& = & (\mu(2\alpha)-\omega')V_\alpha+\omega\int_{c^*\tau}^{+\infty}
\left(\nu \eta_{2\alpha} \zeta_\tau q^2  - \nu (\partial_\tau \eta_{2\alpha}) p q\right)\,dx-2D_\alpha.\no
\end{eqnarray}
Here, as in the case $\omega=0$, we have defined
\[
D_\alpha(\tau)=\int_{c^*\tau}\nu\eta_{2\alpha}\zeta q_x^2dx.
\]
We now use the following fact: for all $M>0$, there is a constant $\kappa_M>0$ such that, for all nonnegative
functions $u(x)\in C^1([0,1])$ such that $\vert u'(x)\vert\leq M\di\int_0^1u(x)dx$ we have:
$$\int_0^1u(x)dx\leq\kappa_M\int_0^1xu(x)dx.$$
If not, there is a sequence $u_n$ of such functions with unit mass and uniformly bounded derivatives whose first moments
tend to 0, an impossibility. Now, from this remark we have
\[
\omega\int_{c^*\tau}^{+\infty} \nu \vert\zeta_\tau\eta_{2\alpha} | q^2 \,dx\leq C\tau^{-1}V_\alpha,
\]
and from Lemma \ref{l6.1} we have
$$
\omega\int_{c^*\tau}^{+\infty}\nu \vert\partial_\tau\eta_{2\alpha} \vert \tilde p q\, dx \leq C \omega\int_{c^*\tau}^{+\infty}\nu \eta_{2\alpha} \tilde p q\, dx \leq C \tau^{-1}V_\alpha.
$$
Because of Lemma \ref{l5.1}, we have (following the lines of the proof of  Proposition~\ref{prop-0523}):
$$
\frac{dV_\alpha}{d\tau}  \leq (\mu(2\alpha)+O(\tau^{-1}))V_\alpha(\tau) - C\frac{V_\alpha^{5/3}}{{I_\alpha}^{4/3}}=
-C\frac{(V_\alpha)^{5/3}}{I_\alpha(0)^{4/3}} e^{\tau R_\alpha} +(\mu(2\alpha)+O(\tau^{-1}))V_\alpha(\tau).
$$
Let us choose $T>0$ and examine the above differential inequality with
$\alpha=T^{-1}$ and $\tau\leq T$. For $\Lambda>0$ large enough, the function $\Lambda\tau^{-3/2}$ is a super-solution for $\tau\leq T$ ,
showing that $V_\alpha(T)=O(T^{-3/2})$. So, for all $\tau>0$, we have $V_\alpha(\tau)\leq C\tau^{-3/2}$,
and the rest of the proof of this proposition follows as in Proposition~\ref{prop-0523}.\qed

%%%%%%%%%%%%%%%%%%%%%%%%%%%%%%%%%%%%%%%%%

\subsection*{The proof of Proposition~\ref{plan-prop2}}

The proof is by a multiple-scale expansion. We will construct a function $\theta^{app}$ having the form
\[
\theta^{app}(\tau,x) = a(\tau)v(\tau,(x - c^*\tau)/R(\tau),x).
\]
which satisfies $\theta^{app}(\tau,c^*\tau)=0$, with $R(\tau) = \tau^{1/2}$. 
Plugging this ansatz into 
$$(1-\omega(\tau))\theta_\tau=\theta_{xx} + 2 \frac{\phi_x}{\phi} \theta_x,$$
we see that $v(\tau,z,x)$ should satisfy
$$
(1-\omega)\left[\frac{a'}{a}v +   v_\tau - \frac{zR'}{R} v_z  - \frac{c^*}{R}v_z\right]
 =\frac{1}{R^2} v_{zz} + \frac{2}{R} v_{zx}  + v_{xx}  
+  2 \frac{\phi_x}{\phi} v_x +  \frac{2}{R}  \frac{\phi_x}{\phi} v_z.
$$
We will construct an approximate solution given by the expansion 
\[
v = v(\tau,z,x) = v^0(z) + \frac{1}{R}v^1(z,x) + \frac{1}{R^2} v^2(z,x)+ \frac{1}{R^3} v^3(z,x),
\]
where $v^1(z,x)$ and $v^2(z,x)$ are uniformly bounded in each compact set in $z$, and $x$, and are both periodic in $x$. Therefore, the desired equality is
\begin{eqnarray}
&  & (1-\omega)\frac{a'}{a}(v^0 + \frac{1}{R} v^1 + \frac{1}{R^2} v^2+\frac{1}{R^3} v^3)  - (1-\omega)\frac{R'}{R^2}\left( v^1 + \frac{2}{R} v^2  + \frac{3}{R^2} v^3\right) \no \\
& & - (1-\omega)\frac{zR'}{R} (v^0_z + \frac{1}{R} v^1_z + \frac{1}{R^2}v^2_z + \frac{1}{R^3}v^3_z)  - (1-\omega)\frac{c^*}{R}(v^0_z + \frac{1}{R} v^1_z+ \frac{1}{R^2}v^2_z+ \frac{1}{R^3}v^3_z) \no \\
&  & \quad \quad \quad = \frac{1}{R^2} v^0_{zz} + \frac{1}{R^3} v^1_{zz} + \frac{1}{R^4} v^2_{zz} + \frac{2}{R^2} v^1_{zx} + \frac{2}{R^3}v^2_{zx}  + \frac{2}{R^4}v^3_{zx}  \no \\ 
& & \quad \quad \quad \quad + \frac{1}{R} v^1_{xx} + \frac{1}{R^2} v^2_{xx} + \frac{1}{R^3} v^3_{xx} + \frac{2 }{R}  \frac{\phi_x}{\phi} v^1_x + \frac{2 }{R^2}  \frac{\phi_x}{\phi} v^2_x + \frac{2 }{R^3}  \frac{\phi_x}{\phi} v^3_x \no \\
& & \quad \quad \quad \quad +  \frac{2}{R}  \frac{\phi_x}{\phi} v^0_z  +  \frac{2}{R^2}  \frac{\phi_x}{\phi} v^1_z+ \frac{2}{R^3}  \frac{\phi_x}{\phi} v^2_z+ \frac{2}{R^4}  \frac{\phi_x}{\phi} v^3_z. \label{ineq1}
\end{eqnarray}
Let us set   $a(\tau) = \tau^{-m}$, so that $a'/a = - m \tau^{-1} = O(R^{-2})$. Now we choose $v_i,\ i\in\{0,...,3\}$ so that terms of order $O(R^{-1})$, $O(R^{-2})$ and $O(R^{-3})$ will cancel. Recall that $\omega(\tau)\sim3/(2c^*\lambda^*\tau)$, so $\omega$ will not play a role until we equate terms of order $O(R^{-3})$, and even then the only term to contribute is $\omega c^*v_z^0/R$. All other terms involving $\omega(\tau)$ are smaller than $O(\tau^{-3/2})$.

If we equate the leading order terms (of order $O(R^{-1})$), we obtain an equation for $v^1$ in terms of $v^0$:
\begin{equation}
 v^1_{xx} + 2  \frac{\phi_x}{\phi} v^1_x = -  \left(  2\frac{\phi_x(x)}{\phi(x)}  + c\right)v^0_z(z).  \label{locancel}
\end{equation}
Recalling $\chi(x)$ defined at (\ref{chidef}) which solves
$$\chi_{xx} + 2 \frac{\phi_x}{\phi} \chi_x = -2\frac{\phi_x}{\phi} - c,$$
we see that (\ref{locancel}) has a solution of the form $v^1(z,x) = v^0_z(z) \chi^0(x) - p^0(z)$ with $\chi^0(x) = \chi(x) + \bar \chi$ being periodic in $x$, and $\bar \chi$ being any constant. For any choice of the constant $\bar \chi$ and $p^0(z)$, (\ref{locancel}) holds and the $O(R^{-1})$ terms in (\ref{ineq1}) cancel.  

Let us now equate the terms of $O(R^{-2})$ in \eqref{ineq1} to obtain
\begin{eqnarray}
 \frac{1}{R^2} v^2_{xx} + \frac{2 }{R^2}  \frac{\phi_x}{\phi} v^2_x  -  \frac{a'}{a}v^0  + \frac{zR'}{R} v^0_z  +\frac{1}{R^2} v^0_{zz}  + c\frac{1}{R^2} v^1_z  + \frac{2}{R^2} v^1_{zx} +  \frac{2}{R^2}  \frac{\phi_x}{\phi} v^1_z = 0, \no
\end{eqnarray}
which is:
\begin{eqnarray}
\label{v2eqn}
  v^2_{xx} + 2  \frac{\phi_x}{\phi} v^2_x  +  m v^0  + \frac{z}{2} v^0_z  + v^0_{zz}  + c v^1_z  + 2 v^1_{zx} +  2  \frac{\phi_x}{\phi} v^1_z = 0.
\end{eqnarray}

Consider the operator $\rho_{xx} + 2 \frac{\phi_x(x)}{\phi(x)} \rho_x = \hat \phi^{-2}(\hat \phi^2 \rho_x)_x$ acting on 1-periodic functions, where $\hat \phi = e^{- \mu x} \psi(x)$. We claim that the adjoint operator has one-dimensional kernel. A function $\eta$ is in the kernel of the adjoint operator if and only if 
\begin{equation}
(\hat \phi^2(\hat \phi^{-2} \eta)_x)_x = 0 \label{adjeta}
\end{equation}
which holds if and only if
\begin{equation}
\eta(x) = k_1 \hat \phi^2(x) \int_0^x \hat \phi^{-2}(s) \,ds + k_2 \hat \phi^2(x) \label{etanull}
\end{equation}
for some constants $k_1$ and $k_2$. If $k_1 = 0$, the function $\eta$ cannot be periodic, since $\hat \phi^2(x) = e^{-2 \mu x} \psi^2(x)$ is not periodic. So, we may assume $k_1 = 1$. However, the function
\[
\eta(x) = \hat \phi^2(x) \int_0^x \hat \phi^{-2}(s) \,ds + k_2 \hat \phi^2(x)
\]
will be periodic for exactly one choice of $k_2$:
\[
k_2 = \frac{\hat \phi^2(1)}{\hat \phi^2(0) - \hat \phi^2(1)} \int_0^1  \hat \phi^{-2}(s) \,ds > 0.
\]
Therefore, with $k_2$ chosen in this way, any other solution of (\ref{adjeta}) must be a multiple of this function $\eta$ given by (\ref{etanull}). Observe that $\eta > 0$ for all $x$. 

If $\eta(x)$ is 1-periodic and spans the kernel of $(\hat \phi^2(\hat \phi^{-2} \eta)_x)_x$, then equation (\ref{v2eqn}) is solvable if and only if the sum
\[
 m v^0  + \frac{z}{2} v^0_z  + v^0_{zz}  + c v^1_z  + 2 v^1_{zx} +  2  \frac{\phi_x}{\phi} v^1_z 
 \]
is orthogonal to $\eta(x)$, for each $z \in \mathbb{R}$. Using $v^1 = v^0_z(z) \chi(x) - p^0(z)$, we write the sum as
\begin{eqnarray}
m v^0  + \frac{z}{2} v^0_z  + v^0_{zz} +  c v^0_{zz}\chi^0 + 2v^0_{zz} \chi^0_x + 2\frac{\phi_x}{\phi} v^0_{zz} \chi^0 - \left( c + 2 \frac{\phi_x}{\phi} \right)p^0_z. \label{rhsum}
\end{eqnarray}
So, the solvability condition is
\begin{eqnarray}
 \left(m v^0  + \frac{z}{2} v^0_z  + v^0_{zz}\right) \int_0^1 \eta(x) \,dx & = &  - \int_0^1 \left(  c v^1_z  + 2 v^1_{zx} +  2  \frac{\phi_x}{\phi} v^1_z  \right) \eta(x) \,dx \no \\
& = & - \int_0^1 \left(  c v^0_{zz}\chi^0 + 2v^0_{zz} \chi^0_x + 2\frac{\phi_x}{\phi} v^0_{zz} \chi^0  \right) \eta(x) \,dx.\no
\end{eqnarray}
Here we have used the fact that $\int_0^1 ( c + 2 \frac{\phi_x}{\phi}) \eta(x) \,dx = 0$, so that the terms involving $p^0_z$ cancel after integration against $\eta$. Hence, $v^0(z)$ should solve
\[
m v^0  + \frac{z}{2} v^0_z  + (1 + \kappa) v^0_{zz} = 0 
\]
where
\be\label{kappadef}
\kappa = \left( \int_0^1 \eta(x) \,dx \right)^{-1} \int_0^1 \left( c \chi^0(x) + 2 \chi^0_x(x) + 2\frac{\phi_x}{\phi} \chi^0(x) \right) \eta(x) \,dx.
\ee
It is not difficult to show that
\[
1 + \kappa = \frac{ \int \eta (1 + \chi^0_x)^2 \,dx}{\int \eta \,dx} > 0
\]
In particular, $\kappa$ is independent of the normalization of $\chi^0(x)$ (the choice of $\bar \chi$). Thus, we choose $v^0(z) > 0$ to be the principal eigenfunction of
$$m v^0  + \frac{z}{2} v^0_z  + (1 + \kappa) v^0_{zz} =0, \quad z>0,\quad\quad\quad v^0(0)=0,$$
which forces $m=1$, and
\[
v^0(z)=ze^{-\frac{z^2}{4(1 + \kappa)}}.
\]

The function $p^0(z)$ is undetermined so far.  With $v^0(z)$ chosen in this way, there exists a function $v^2(z,x)$ which is periodic in $x$ and satisfies (\ref{v2eqn}). Thus, the $O(R^{-2})$ terms cancel. In consideration of (\ref{rhsum}) and the definition of $v^0$, we see that (\ref{v2eqn}) is equivalent to
\[
v^2_{xx} + 2  \frac{\phi_x}{\phi} v^2_x   =  -v^0_{zz}(z) \left(  c  \chi^0 + 2 \chi^0_x + 2\frac{\phi_x}{\phi} \chi^0  - \kappa \right)- \left( c + 2 \frac{\phi_x}{\phi} \right)p^0_z.
\]
Therefore, $v^2(z,x)$ must have the form
\[
v^2(z,x) = v^0_{zz}(z) \hat v^2(x) - p^0_z(z) \chi^0(x) + p^1(z),
\]
where $\hat v^2(x)$ is a periodic solution of
\[
\hat v^2_{xx} + 2  \frac{\phi_x}{\phi} \hat v^2_x = - \left(  c  \chi^0 + 2 \chi^0_x + 2\frac{\phi_x}{\phi} \chi^0  - \kappa \right).
\]
Finally, equating the $R^{-3}$ terms suggests choosing $v^3(x,z)$ to satisfy
\begin{equation}
v^3_{xx}+2\frac{\phi_x}{\phi}v^3_x=\frac3{2\lambda^*}v^0_z-(m+1)v^1-\frac{z}2v^1_z-v^1_{zz}-(c^*+2\frac{\phi_x}{\phi})v_z^2-2v^2_{zx}. \label{v3eqn}
\end{equation}
The right hand side is:
\begin{eqnarray}
&&\frac3{2\lambda^*}v^0_z - 2 v^0_z \chi^0 + 2 p^0 - \frac{z}{2} v^0_{zz} \chi^0 + \frac{z}{2} p^0_{z} - v^0_{zzz} \chi^0 + p^0_{zz} \no \\
& & \quad \quad -(c^*+2\frac{\phi_x}{\phi}) \left( v^0_{zzz} \hat v^2 - p^0_{zz} \chi^0 + p^1_z \right) - 2 \left( v^0_{zzz} \hat v^2_x - p^0_{zz} \chi^0_x \right)\no\end{eqnarray}
Therefore, the solvability condition implies that $p^0(z)$ should satisfy
$$2 p^0  + \frac{z}{2} p^0_{z}  + (1 + \kappa) p^0_{zz}  = \beta_1 v^0_{zzz} + \beta_2 \frac{z}{2} v^0_{zz} + (\frac3{2\lambda^*} - 2 \beta_2) v^0_z$$
where
\begin{eqnarray*}
&&\beta_1 = \left ( \int_0^1 \eta(x) \,dx \right)^{-1} \int_0^1 \left(\chi^0 + (c^*+2\frac{\phi_x}{\phi})\hat v^2 + 2 \hat v^2_x \right)\eta(x)  \,dx, 
\\
&&
\beta_2 = \left ( \int_0^1 \eta(x) \,dx \right)^{-1}  \int_0^1 \chi^0 \eta \,dx , 
\end{eqnarray*}
and we would like to have $p^0(0) = 0$. The $p^1$ term does not appear in the solvability condition. 
Therefore, we may take $p^1(z) \equiv 0$. We let $p^0(z)$ be the unique solution of the initial value problem
$$2 p^0  + \frac{z}{2} p^0_{z}  + (1 + \kappa) p^0_{zz}  = \beta_1 v^0_{zzz} + \beta_2 \frac{z}{2} v^0_{zz} + 
(\frac3{2\lambda^*} - 2 \beta_2) v^0_z, \quad z > 0$$
with the initial data initial $p^0(z) = 0$ and $p^0_z(0) = 0$. 

Having chosen $p^0$ in this way, we take $v^3$ to be a solution of (\ref{v3eqn}), which is unique up to addition of a 
function $p^3(z)$. So, the $O(R^{-3}) = O(\tau^{-3/2})$ terms have canceled. Our approximate solution is:
\[
\theta^{app}(t,x)  =  \tau^{-1} v^0(z) + \tau^{-3/2} v^1(z,x) + \tau^{-2} v^2(z,x)  + \tau^{-5/2} v^3(z,x),
\]
with
\begin{eqnarray}
v^0(z) & = & z e^{-\frac{z^2}{4(1 + \kappa)}},\no \\
v^1(z,x) & = & \chi^0(x)e^{-\frac{z^2}{4(1 + \kappa)}} - \frac{z^2\chi^0(x)}{2(1 + \kappa)}e^{-\frac{z^2}{4(1 + \kappa)}} - p^0(z). \no 
\end{eqnarray}
Now, fix a constant $\sigma > 0$. Having chosen $p^0(0) = 0$ and $p^0_z(0) = 0$, we may choose $C_1 > 0$ so that $|p^0(z)| \leq C_1 z^2$ for all $z \in [0,\sigma]$. Consequently, there is a constant $C_2 > 0$ such that for all $x \in [c^* \tau, c^* \tau + \sigma \sqrt{\tau}]$ and $\tau > 1$ we have 
\[
\left| \theta^{app}(t,x) - \frac{x - c^* \tau + \chi^0(x)}{\tau^{3/2}} e^{-\frac{(x - c^* \tau)^2}{4(1 + \kappa)\tau}} \right| \leq C_2 \tau^{-3/2} \left( \frac{x - c^* \tau}{\sqrt{\tau}} \right)^2 + O(\tau^{-2})
\]
The last term $O(\tau^{-2})$ comes from $v^2$ and $v^3$ and the fact that $v^2$ and $v^3$ are uniformly bounded over $(z,x) \in [0,\sigma] \times \Rm$. 

Since the periodic function $\chi^0(x) = \chi(x) + \bar \chi$ is unique up to addition of a constant, we may choose $\bar \chi < 0$ so that $\max_{x} \chi^0(x) < -1$. Then, at the point $x = c^* \tau$ we have
\[
 \theta^{app}(t,c^* \tau) \leq \tau^{-3/2} \chi^0(c^* \tau) + O(\tau^{-2}) \leq - \tau^{-3/2} + O(\tau^{-2}), 
\]
which is negative for all $\tau > 1$ sufficiently large.  Alternatively, we could choose $\bar \chi > 0$ so that $\min_x \chi^0(x) > 0$. Then we would have $\theta^{app}(\tau, c^* \tau) > 0$ for all $\tau$ sufficiently large. \qed

%%%%%%%%%%%%%%%%%%%%%%%%%%%%%%%%%%%%%%%%%

\subsection*{The proof of Proposition~\ref{plan-prop3}}

Using Lemma~\ref{lem:selfadj} we bring this problem into the form
\begin{equation}\label{plan-eq42}
(1-\omega(\tau))\xi_\tau=\frac{1}{\nu(x)}\pdr{}{x}\left(\nu(x) \xi_x\right)-\frac{c^*}{\nu(x)}\xi_x.
\end{equation}
Let
\[
\Phi(\tau,x)=\xi(\tau,x)-\theta^{app}(\tau,x)
\]
so that $\Phi(\tau,c^* \tau) = 0$ and $\Phi(\tau,c^* \tau + L_0 + \varepsilon \sqrt{\tau}) = 0$. We have
\[
(1-\omega(\tau))\nu(x)\Phi_\tau=(\nu(x)\Phi_x)_x-c^*\Phi_x+ O(\tau^{-3}).
\]
Multiplying by $\Phi(\tau,x)$ and integrating by parts over the interval $I=[c^*\tau,c^*\tau+ L_0 + \eps\sqrt{\tau}]$, we obtain
\[
\farc{1}{2}\frac{d}{d\tau}\int_I\nu(x)(1-\omega(\tau))\Phi^2dx+\farc{\omega'(\tau)}{2}\int_I\nu(x)\Phi^2 dx=
-\int_{I}\nu(x)\Phi_x^2dx+O(\tau^{-3})\int_I\Phi dx.
\]
Note that, since $\Phi(\tau,c^*\tau)=0$, we have
\[
{|\omega'(\tau)|}\int_I\nu\Phi^2dx\le \farc{C}{\tau^2}\eps^2\tau\int_I\nu\Phi_x^2dx
\]
and
$$
{|O(\tau^{-3})\int_I\Phi dx|\le \farc{C\varepsilon}{\tau^{9/2}}}+\frac1\tau\int_I\Phi^2dx\le \frac{C}{\tau^4} + C\eps^2\int_I\nu\Phi_x^2dx.
$$
If now $\eps$ is small enough so that the constant $C \eps^2$ is less than $1/4$ it follows  that, for $\tau>\tau_0$ large enough, we have
\begin{eqnarray}
\farc{1}{2}\frac{d}{d\tau}\int_I\nu(x)(1-\omega(\tau))\Phi^2dx & \le &
-\farc{1}{2}\int_{I}\nu(x)\Phi_x^2dx+\frac{C}{\tau^{4}}\no\\
& \le & -\farc{1}{C(L_0+\eps\sqrt{\tau})^2}\int_{I}\nu(x)
(1-\omega(\tau))\Phi^2dx+\frac{C}{\tau^{4}}.\no
\end{eqnarray}
We conclude that, for $\eps$ sufficiently small, we have
\[
\int_I\nu(x)\Phi^2dx\le\frac{C_\eps}{(1+\tau)^{1/\eps^2}}+\frac{C_\eps}{(1+\tau)^{3}}.
\]
Now, parabolic regularity implies that $|\Phi(\tau,x)|\le C/(1+\tau)^{3/2}$ for $\tau>\tau_0$ sufficiently large. This completes
the proof of Proposition~\ref{plan-prop3}.~$\Box$

%%%%%%%%%%%%%%%%%%%%%%%%%%%%%%%%%%%%%%%%%
%%%%%%%%%%%%%%%%%%%%%%%%%%%%%%%%%%%%%%%%%

\section{Proof of Theorem~\ref{th2}}

This section is devoted to the proof of the convergence of the solution $u$ to the family of shifted minimal fronts $U_{c^*}$. We first remember that $u$ is bounded away from $0$ or $1$ around the position $c^*t-(3/(2\lambda^*))\ln t$ for large $t$. To the right of this position, the solution $u$ has the same type of decay as the critical front $U_{c^*}$, as it follows from the estimates of Sections~2 and~3. Therefore,~$u$ is almost trapped between two finite shifts of the profile of the front $U_{c^*}$. From a Liouville-type result, similar to that in~\cite{bh4} and based on the sliding method, the convergence to the shifted approximated minimal fronts will follow.

First, we derive from Sections~2 and~3 some exponential bounds of $u$ to the right of the position $c^*t-(3/(2\lambda^*))\log t$.

\begin{lem}\label{lemkapparho} Let $\sigma>0$ be as in Proposition~$\ref{cor-0524}$. There exist two positive constants $0<\kappa\le\rho$ such that
\be\label{lowerkappa}
\kappa\,y\,e^{-\lambda^*y}\le u\big(t,c^*t-\frac{3}{2\lambda^*}\log t+y\big)\ \hbox{ for all }t\ge 1\hbox{ and }0\le y\le\sigma\sqrt{t}
\ee
and
\be\label{upperrho}
u\big(t,c^*t-\frac{3}{2\lambda^*}\log t+y\big)\le\rho\,y\,e^{-\lambda^*y}\ \hbox{ for all }t\ge 1\hbox{ and }y\ge 1.
\ee
\end{lem}

\noindent{\bf{Proof.}} Under the notations of Section~2, it follows from~(\ref{deftildeU}),~(\ref{080414}) and~(\ref{defr}) that, for some positive constants $T_1$ and $L_0$,
\be\label{lower1}
u(t,x)\ge U_{c^*}^k\big(t-\frac{3}{2\lambda^*c^*}\log t-L_0,x\big)\ \hbox{ for all }t\ge T_1\hbox{ and }0\le x\le c^*t+\sigma\sqrt{t}.
\ee
The pulsating front $U_{c^*}^k$ can be written as $U_{c^*}^k(t,x)=\phi_{c^*}^k(x-c^*t,x)$, where $0<\phi_{c^*}^k(s,x)<1$ is continuous in $\R\times\R$, $1$-periodic in $x$, and $\phi_{c^*}^k(-\infty,\cdot)=1$, $\phi_{c^*}^k(+\infty,\cdot)=0$. Furthermore, it is known~\cite{Hamel1} that there is a constant $\beta>0$ such that $\phi_{c^*}^k(s,x)\sim\beta\,\psi(x,\lambda^*)\,s\,e^{-\lambda^*s}$ as $s\to+\infty$, uniformly in $x\in\R$. In particular, there is $\eta>0$ such that $\phi_{c^*}^k(s,x)\ge\max\big(\eta\,s\,e^{-\lambda^*s},0\big)$ for all $(s,x)\in\R\times\R$. As a consequence, it follows from~(\ref{lower1}) that
$$u\big(t,c^*t-\frac{3}{2\lambda^*}\log t+y\big)\ge\phi_{c^*}^k(y+c^*L_0,c^*t-\frac{3}{2\lambda^*}\log t+y\big)\ge\eta\,(y+c^*L_0)\,e^{-\lambda^*(y+c^*L_0)}$$
for all $t\ge T_1$ and $-c^*t+(3/(2\lambda^*))\log t\le y\le\sigma\sqrt{t}+(3/(2\lambda^*))\log t$. Therefore, there are $T_2\ge T_1$ and $\kappa>0$ such that
$$u\big(t,c^*t-\frac{3}{2\lambda^*}\log t+y\big)\ge\kappa\,y\,e^{-\lambda^*y}\ \hbox{ for all }t\ge T_2\hbox{ and }0\le y\le\sigma\sqrt{t}.$$
The inequality~(\ref{lowerkappa}) follows, by positivity and continuity of $u$ over $[1,+\infty)\times\R$, by taking a smaller $\kappa>0$ if necessary.\par
On the other hand, it follows from~(\ref{plan-eq14}),~(\ref{defbaru}) and~(\ref{ubaru}) that there exist some positive constants $\overline{T}$, $\overline{y}$ and $\rho$ such that
$$u\big(t,c^*t-\frac{3}{2\lambda^*}\log t+y\big)\le\rho\,y\,e^{-\lambda^*y}\ \hbox{ for all }t\ge\overline{T}\hbox{ and }y\ge\overline{y}.$$
The inequality~(\ref{upperrho}) follows, by positivity and continuity of $u$ over $[1,+\infty)\times\R$, by taking a larger $\rho>0$ if necessary.~\qed\hfill\break

\noindent{\bf{Proof of Theorem~\ref{th2}.}} First, let $\sigma>0$ and $0<\kappa\le\rho$ be given as in the previous lemma. Write the pulsating front $U_{c^*}$ as
\be\label{defphic*}
U_{c^*}(t,x)=\phi_{c^*}(x-c^*t,x),
\ee
where $0<\phi_{c^*}(s,x)<1$ is continuous in $\R\times\R$, $1$-periodic in $x$, and $\phi_{c^*}(-\infty,\cdot)=1$, $\phi_{c^*}(+\infty,\cdot)=0$. From~\cite{Hamel1}, there is a constant $B>0$ such that
\be\label{Uc*}
\phi_{c^*}(s,x)\sim B\,\psi(x,\lambda^*)\,s\,e^{-\lambda^*s}\ \hbox{ as }s\to+\infty,\hbox{ uniformly in }x\in\R.
\ee
Choose now any real number $\widetilde{C}\ge 0$ so that
\be\label{defC}
B\,\max\psi(\cdot,\lambda^*)\,e^{-c^*\lambda^*\widetilde{C}}\le\kappa\le\rho\,e^{\lambda^*}\le B\,\min\psi(\cdot,\lambda^*)\,e^{c^*\lambda^*\widetilde{C}}.
\ee\par
Let us prove that~(\ref{fronts}) holds with the choice of $C=\widetilde{C}+1/c^*$. Assume not. There are then $\varepsilon>0$ and a sequence of positive times $(t_n)_{n\in\N}$ such that $t_n\to+\infty$ as $n\to+\infty$ and
\[
\min_{|\xi| \leq\widetilde{C}+1/c^*} \left \Vert u(t_n,\cdot) - U_{c^*}\big(t_n-\frac{3}{2c^*\lambda^*}\log t_n+\xi,\cdot\big) \right \Vert_{L^\infty(0,+\infty)} \ge\varepsilon
\]
for all $n\in\N$. Since $\phi_{c^*}(-\infty,\cdot)=1$, $\phi_{c^*}(+\infty,\cdot)=0$ uniformly in $\R$ and $\phi(s,x)$ is $1$-periodic in $x$, it follows from~(\ref{defphic*}) and Theorem~\ref{thm-delay-apr5} that there exists a constant $\theta\ge0$ such that
\be\label{minmax}
\min_{|\xi|\le\widetilde{C}}\Big(\max_{|y|\le\theta}\Big|u\big(t_n,y+\big[c^*t_n-\frac{3}{2\lambda^*}\ln t_n\big]\big)-U_{c^*}(\xi,y)\Big|\Big)\ge\varepsilon
\ee
for all $n\in\N$, where $[c^*t_n-3/(2\lambda^*)\log t_n]$ denotes the integer part of $c^*t_n-3/(2\lambda^*)\log t_n$.\par
For each $n\in\N$, set
$$u_n(t,x)=u\big(t+t_n,x+\big[c^*t_n-\frac{3}{2\lambda^*}\log t_n\big]\big).$$
Up to extraction of a subsequence, the functions $u_n$ converge locally uniformly in $\R^2$ to a solution $u_{\infty}$ of
\be\label{uinfty}
(u_{\infty})_t=(u_{\infty})_{xx}+g(x)\,f(u_{\infty})\ \hbox{ in }\R^2
\ee
such that $0\le u_{\infty}\le 1$ in $\R^2$. Furthermore, Theorem~\ref{thm-delay-apr5} implies that
$$\lim_{A\to+\infty}\Big(\sup_{(t,x)\in\R^2,\,x \geq c^*t+A}u_{\infty}(t,x)\Big)=0$$
and
\be\label{lim1}
\lim_{A\to-\infty}\Big(\inf_{(t,x)\in\R^2,\,x\le c^*t + A}u_{\infty}(t,x)\Big)=1.
\ee\par
On the other hand, for each fixed $t\in\R$ and $y>2$, and $n$ large enough, write
$$u_n(t,c^*t+y)=u\big(t+t_n,c^*(t+t_n)-\frac{3}{2\lambda^*}\log(t+t_n)+y+\gamma_n\big),$$
where
$$\gamma_n=\big[c^*t_n-\frac{3}{2\lambda^*}\log t_n\big]-\big(c^*t_n-\frac{3}{2\lambda^*}\log(t+t_n)\big).$$
There holds $t+t_n\ge 1$ and $1\le y+\gamma_n\le\sigma\sqrt{t+t_n}$ for $n$ large enough, whence
$$\kappa\,(y+\gamma_n)\,e^{-\lambda^*(y+\gamma_n)}\le u_n(t,c^*t+y)\le\rho\,(y+\gamma_n)\,e^{-\lambda^*(y+\gamma_n)}$$
for $n$ large enough, from Lemma~\ref{lemkapparho}. Since $-1\le\liminf_{n\to+\infty}\gamma_n\le\limsup_{n\to+\infty}\gamma_n\le 0$, it follows that
\be\label{exp}
\kappa\,(y-1)\,e^{-\lambda^*y}\le u_{\infty}(t,c^*t+y)\le\rho\,y\,e^{-\lambda^*(y-1)}\ \hbox{ for all }t\in\R\hbox{ and }y\ge2.
\ee\par
The following Liouville-type result gives a classification of the time-global solutions $u_{\infty}$ of~(\ref{uinfty}) satisfying the above properties~(\ref{lim1}) and~(\ref{exp}).

\begin{lem}\label{lemliouville}
For any solution $0\le u_{\infty}\le 1$ of~$(\ref{uinfty})$ in $\R^2$ satisfying~$(\ref{lim1})$ and~$(\ref{exp})$ for some positive constants $\kappa$ and $\rho$, there is $\xi_0\in\R$ such that
\be\label{uinftybis}
u_{\infty}(t,x)=U_{c^*}(t+\xi_0,x)\ \hbox{ for all }(t,x)\in\R^2.
\ee
\end{lem}

The proof of this lemma is postponed at the end of this section. We first complete the proof of Theorem~\ref{th2}. It follows from Lemma~\ref{lemliouville}, from~(\ref{defphic*}), from~(\ref{exp}) and from the exponential decay~(\ref{Uc*}) of $\phi_{c^*}$, that
$$\kappa\le B\max\psi(\cdot,\lambda^*)\,e^{c^*\lambda^*\xi_0}\ \hbox{ and }\ B\min\psi(\cdot,\lambda^*)\,e^{c^*\lambda^*\xi_0}\le\rho\,e^{\lambda^*},$$
whence $|\xi_0|\le\widetilde{C}$ from~(\ref{defC}). But since (at least for a subsequence) $u_n\to u_{\infty}$ locally uniformly in $\R^2$, it follows in particular that $u_n(0,\cdot)-U_{c^*}(\xi_0,\cdot)\to 0$ uniformly in $[-\theta,\theta]$, that is
$$\max_{|y|\le\theta}\Big|u\big(t_n,y+\big[c^*t_n-\frac{3}{2\lambda^*}\log t_n\big]\big)-U_{c^*}(\xi_0,y)\Big|\to0\ \hbox{ as }n\to+\infty.$$
Since $|\xi_0|\le\widetilde{C}$, one gets a contradiction with~(\ref{minmax}). Therefore,~(\ref{fronts}) is proved.\par
Let us now turn to the proof of~(\ref{fronts2}). Let $m\in(0,1)$ be fixed and let $(t_n)_{n\in\N}$ and $(x_n)_{n\in\N}$ be two sequences of positive real numbers such that $t_n\to+\infty$ as $n\to+\infty$ and $u(t_n,x_n)=m$ for all $n\in\N$. Set
$$X_n=[x_n]-\big[c^*t_n-\frac{3}{2\lambda^*}\log t_n\big].$$
Theorem~\ref{thm-delay-apr5} implies that the sequence of integers $(X_n)_{n\in\N}$ is bounded, and may then be assumed to be equal to a constant integer $X_{\infty}$, up to extraction of a subsequence. Under the notations of the previous paragraphs, the functions
$$v_n(t,x)=u(t+t_n,x+[x_n])=u\big(t+t_n,x+X_{\infty}+\big[c^*t_n-\frac{3}{2\lambda^*}\log t_n\big]\big)=u_n(t,x+X_{\infty})$$
converge locally uniformly in $\R^2$, up to extraction of another subsequence, to the function
$$v_{\infty}(t,x)=u_{\infty}(t,x+X_{\infty})=U_{c^*}(t+\xi,x+X_{\infty})=U_{c^*}\Big(t+\xi-\frac{X_{\infty}}{c^*},x\Big)$$
for some real number $\xi$. Since $v_n(0,x_n-[x_n])=m$ for all $n\in\N$ and $x_n-[x_n]\to x_{\infty}$ as $n\to+\infty$, one gets that $U_{c^*}(\xi-X_{\infty}/c^*,x_{\infty})=m$, that is $\xi-X_{\infty}/c^*=T$, where $T$ is the unique real number such that $U_{c^*}(T,x_{\infty})=m$. Finally, the limit $v_{\infty}$ is uniquely determined and the whole sequence~$(v_n)_{n\in\N}$ therefore converges to the pulsating front $U_{c^*}(t+T,x)$. The proof of Theorem~\ref{th2} is thereby complete.~$\Box$\hfill\break

\noindent{\bf{Proof of Lemma~\ref{lemliouville}.}} In the homogeneous case, if, instead of~(\ref{lim1}) and~(\ref{exp}), the function $u_{\infty}$ is assumed to be trapped between two shifts of the minimal traveling front, then the conclusion follows directly from Theorem~3.5 of~\cite{bh4}. In our periodic case, the comparisons~(\ref{exp}) and the exponential behavior~(\ref{Uc*}) of the minimal front~$U_{c^*}$ imply that $u_{\infty}$ is actually trapped between two finite time-shifts of~$U_{c^*}$ in the region $\big\{x-c^*t\ge0\big\}$. In the region where $x-c^*t$ is very negative,~$u_{\infty}$ is close to~$1$ and the maximum principle holds, from the negativity of $f'(1)$: the solution $u_{\infty}$ can then be compared to some of its shifts in this region. We finally complete the proof of the lemma by using a sliding method: we shift the function~$u_{\infty}(t,x+1)$ in time, we compare it with the function $u_{\infty}$, and we show that $u_{\infty}(t+1/c^*,x+1)=u_{\infty}(t,x)$ in $\R^2$. Together with~(\ref{lim1}) and~(\ref{exp}), this will mean that $u_{\infty}$ is a pulsating front. From the uniqueness of the pulsating fronts up to time-shifts~\cite{hr}, the conclusion~(\ref{uinftybis}) will follow. More precisely, for all $\xi\in\R$ and $(t,x)\in\R^2$, we set
$$v^{\xi}(t,x)=u_{\infty}(t+\xi,x+1).$$
We shall compare $v^{\xi}$ to $u_{\infty}$ and prove that $v^{\xi}\ge u_{\infty}$ in $\R^2$ for all $\xi$ large enough. We will then prove that $v^{\xi}\equiv u_{\infty}$ in $\R^2$ for the smallest such $\xi$, and finally that this critical shift is equal to $1/c^*$.\par
To do so, we first notice that, for all $a\le b\in\R$, there holds
\be\label{infsup}
0<\inf_{(t,x)\in\R^2,\,a\le x-c^*t\le b}u_{\infty}(t,x)\le\sup_{(t,x)\in\R^2,\,a\le x-c^*t\le b}u_{\infty}(t,x)<1.
\ee
This a consequence of the strong maximum principle, parabolic regularity, and the fact the solution $0<u_{\infty}<1$ converges to two different limits ($0$ and $1$) as $x-c^*t\to\pm\infty$. Let now $\delta\in(0,1)$ be such that $f$ is nonincreasing in $[1-\delta,1]$, and let us extend $f$ by $0$ on $(1,+\infty)$. From~(\ref{lim1}), there is $A>0$ such that
\be\label{delta}
u_{\infty}(t,x)\ge 1-\delta\hbox{ for all }(t,x)\in\R^2\hbox{ such that }x-c^*t\le-A.
\ee\par
As far as the region $\big\{x-c^*t\ge-A\big\}$ is concerned, we claim that there is $\overline{\xi}\in\R$ such that
\be\label{xi0}
v^{\xi}(t,x)\ge u_{\infty}(t,x)\ \hbox{ for all }x-c^*t\ge-A\hbox{ and }\xi\ge\overline{\xi}.
\ee
Assume not. Then there exist some sequences $(\xi_n)_{n\in\N}$ in $[0,+\infty)$ and $(t_n,x_n)_{n\in\N}$ in $\R^2$ such that $\lim_{n\to+\infty}\xi_n=+\infty$ and
$$x_n-c^*t_n\ge-A,\ \ u_{\infty}(t_n+\xi_n,x_n+1)=v^{\xi_n}(t_n,x_n)<u_{\infty}(t_n,x_n)\ \hbox{ for all }n\in\N.$$
Because of~(\ref{lim1}),~(\ref{exp}) and~(\ref{infsup}), the sequence $(x_n-c^*t_n-c^*\xi_n)_{n\in\N}$ is bounded from below by a constant $M$. Thus,~(\ref{exp}) and~(\ref{infsup}) provide the existence of some positive constants $\widetilde{\kappa}$ and $\widetilde{\rho}$ such that
\be\label{xin}\baa{rcll}
\widetilde{\kappa}\,(x_n-c^*t_n-c^*\xi_n-M+1)\,e^{-\lambda^*(x_n-c^*t_n-c^*\xi_n)} & \le & u_{\infty}(t_n+\xi_n,x_n+1)\vspace{3pt}\\
& < & u_{\infty}(t_n,x_n)\vspace{3pt}\\
& \le & \widetilde{\rho}\,(x_n-c^*t_n+A+1)\,e^{-\lambda^*(x_n-c^*t_n)}\eaa
\ee
for all $n\in\N$. On the other hand, 
$$\baa{rcl}
x_n-c^*t_n+A+1 & = & (x_n-c^*t_n-c^*\xi_n-M+1)+(c^*\xi_n+M+A)\vspace{3pt}\\
& \le & 2\,(x_n-c^*t_n-c^*\xi_n-M+1)\,(c^*\xi_n+M+A)\eaa$$
for $n$ large enough. Putting this into~(\ref{xin}) and passing to the limit as $n\to+\infty$ (with $\xi_n\to+\infty$ as $n\to+\infty$) leads to a contradiction. Thus, the claim~(\ref{xi0}) is proved.\par
Without loss of generality, one can assume that $\overline{\xi}\ge 1/c^*$. In this paragraph, we fix $\xi$ in the interval $[\overline{\xi},+\infty)$. Set
$$\varepsilon^*=\min\Big\{\varepsilon\ge 0,\ v^{\xi}(t,x)+\varepsilon\ge u_{\infty}(t,x)\hbox{ for all }(t,x)\in\R^2\hbox{ such that }x-c^*t\le -A\Big\}$$
and let us prove that $\varepsilon^*=0$. Assume that $\varepsilon^*>0$. Since $u_{\infty}$ is globally Lipschitz continuous and since $v^{\xi}\ge u_{\infty}$ on $\big\{x-c^*t=-A\big\}$ and both functions $v^{\xi}$ and $u_{\infty}$ converge to $1$ as $x-c^*t\to-\infty$, there are a sequence of positive real numbers $(\varepsilon_n)_{n\in\N}$, a sequence $(t_n,x_n)_{n\in\N}$ in $\R^2$ and a real number $y_{\infty}<-A$ such that
$$\varepsilon_n\to\varepsilon^*,\ x_n-c^*t_n\to y_{\infty}\hbox{ as }n\to+\infty\hbox{ and }v^{\xi}(t_n,x_n)+\varepsilon_n<u_{\infty}(t_n,x_n)\hbox{ for all }n\in\N.$$
Without loss of generality, one can also assume that
$$x_n-[x_n]\to x_{\infty}\hbox{ and }t_n-\frac{[x_n]}{c^*}\to\tau\hbox{ as }n\to+\infty,$$
with $y_{\infty}=x_{\infty}-c^*\tau$. Up to extraction of a subsequence, the functions
$$U_n(t,x)=u_{\infty}\Big(t+\frac{[x_n]}{c^*},x+[x_n]\Big)$$
converge locally uniformly in $\R^2$ to a solution $U_{\infty}$ of~(\ref{uinfty}) satisfying~(\ref{lim1}) and~(\ref{exp}). Set
$$V^{\xi}(t,x)=U_{\infty}(t+\xi,x+1)\hbox{ for all }(t,x)\in\R^2.$$
Therefore, $V^{\xi}(t,x)+\varepsilon^*\ge U_{\infty}(t,x)$ for all $(t,x)\in\R^2$ such that $x-c^*t\le-A$, with equality at the point $(\tau,x_{\infty})$ such that $x_{\infty}-c^*\tau=y_{\infty}<-A$. On the other hand, for all $(t,x)\in\R^2$ such that $x-c^*t\le -A$, there holds
$$V^{\xi}(t,x)+\varepsilon^*\ge V^{\xi}(t,x)\ge 1-\delta$$
from~(\ref{delta}), the definition of the functions $V^{\xi}$ and $U_n$, and the assumption $\xi\ge1/c^*$. Consequently,
$$V^{\xi}_t(t,x)-V^{\xi}_{xx}(t,x)=g(x)f(V^{\xi}(t,x))\ge g(x)f(V^{\xi}(t,x)+\varepsilon^*)$$
for all $(t,x)\in\R^2$ such that $x-c^*t\le-A$, since $f$ is nonincreasing in $[1-\delta,+\infty)$ and $g$ is positive. Since $U_{\infty}$ solves~(\ref{uinfty}), it follows from the strong maximum principle that $V^{\xi}(t,x)+\varepsilon^*=U_{\infty}(t,x)$ for all $(t,x)\in\R^2$ such that $x-c^*t\le-A$ and $t\le\tau$. The positivity of $\varepsilon^*$ is in contradiction with the fact that $V^{\xi}$ and $U_{\infty}$ converge to $1$ uniformly as $x-c^*t\to-\infty$. Therefore, $\varepsilon^*=0$, whence
\be\label{epsilon*}
v^{\xi}(t,x)\ge u_{\infty}(t,x)\hbox{ for all }(t,x)\in\R^2\hbox{ such that }x-c^*t\le-A.
\ee
Together with~(\ref{xi0}), one gets finally that $v^{\xi}\ge u_{\infty}$ in $\R^2$ for all $\xi\ge\overline{\xi}$.\par
Set now
$$\xi_*=\min\Big\{\xi\in\R,\ v^{\xi'}\ge u_{\infty}\hbox{ in }\R^2\hbox{ for all }\xi'\ge\xi\Big\},$$
which is a well defined real number such that $\xi_*\le\overline{\xi}$ (notice that $v^{\xi}(t,x)\to0$ as $\xi\to-\infty$ for each fixed $(t,x)\in\R^2$, while $u_{\infty}>0$ in $\R^2$). Our goal is to prove that
$$\xi_*\le\frac{1}{c^*},$$
which will then yield $v^{1/c^*}\ge u_{\infty}$ and a symmetric argument will then give the desired conclusion.\par
Assume then that $\xi_*>1/c^*$. Remember that $v^{\xi_*}\ge u_{\infty}$ by definition of $\xi_*$. We first claim that, for any $a\le b$ in $\R$,
\be\label{infab}
\inf_{(t,x)\in\R^2,\,a\le x-c^*t\le b}\big(v^{\xi_*}(t,x)-u_{\infty}(t,x)\big)>0.
\ee
Otherwise, by a usual limiting argument, there would exist a solution $0\le U_{\infty}\le 1$ of~(\ref{uinfty}) satisfying~(\ref{lim1}) and~(\ref{exp}), and such that $U_{\infty}(t+\xi_*,x+1)\ge U_{\infty}(t,x)$ for all $(t,x)\in\R^2$ with equality somewhere. From the strong maximum principle and the uniqueness of the solutions of the Cauchy problem associated to~(\ref{uinfty}), it would then follow that $U_{\infty}(t+\xi_*,x+1)=U_{\infty}(t,x)$ for all $(t,x)\in\R^2$ and then $U_{\infty}(t+k\xi_*,x+k)=U_{\infty}(t,x)$ in $\R^2$ for all $k\in\N$. Since one has assumed that $\xi_*>1/c^*$ and since $U_{\infty}$ satisfies~(\ref{lim1}), the limit as $k\to+\infty$ implies that $U_{\infty}(t,x)=1$ for all $(t,x)\in\R^2$, which is clearly impossible, because of property~(\ref{exp}) satisfied by~$U_{\infty}$.\par
Therefore,~(\ref{infab}) holds. In particular, since $u_{\infty}$ is Lipschitz, there is $\underline{\xi}\in(1/c^*,\xi_*)$ such that
$$v^{\xi}(t,x)\ge u_{\infty}(t,x)\hbox{ for all }(t,x)\in\R^2\hbox{ such that }x-c^*t=-A\hbox{ and for all }\xi\in[\underline{\xi},\xi_*].$$
Furthermore, $v^{\xi}(t,x)\ge1-\delta$ for all $(t,x)\in\R^2$ such that $x-c^*t\le-A$ and for all $\xi\in[\underline{\xi},\xi_*]\subset[1/c^*,+\infty)$, from~(\ref{delta}) and the definition of $v^{\xi}$. As done in the proof of~(\ref{epsilon*}), it follows then that
\be\label{xi}
v^{\xi}(t,x)\ge u_{\infty}(t,x)\hbox{ for all }(t,x)\in\R^2\hbox{ such that }x-c^*t\le-A\hbox{ and for all }\xi\in[\underline{\xi},\xi_*].
\ee\par
On the other hand, the definition of $\xi_*$ implies that there exist a sequence $(\xi_n)_{n\in\N}$ in $(\xi_*-1,\xi_*)$ and a sequence $(t_n,x_n)_{n\in\N}$ in $\R^2$ such that
\be\label{xinbis}
\xi_n\to\xi_*\hbox{ as }n\to+\infty\hbox{ and }v^{\xi_n}(t_n,x_n)<u_{\infty}(t_n,x_n)\hbox{ for all }n\in\N.
\ee
Property~(\ref{xi}) yields $x_n-c^*t_n>-A$ for all $n$ large enough and~(\ref{infab}) and~(\ref{xinbis}) imply then that $x_n-c^*t_n\to+\infty$ as $n\to+\infty$. Up to extraction of a subsequence, one can assume that $x_n-[x_n]\to x_{\infty}\in[0,1]$ as $n\to+\infty$.\par
Define now
$$U_n(t,x)=\frac{u_{\infty}(t+t_n,x+[x_n])}{u_{\infty}(t_n,[x_n])}\hbox{ and }V_n(t,x)=\frac{v^{\xi_*}(t+t_n,x+[x_n])}{u_{\infty}(t+t_n,x+[x_n])}$$
for all $(t,x)\in\R^2$ and $n\in\N$. From~(\ref{exp}) and $\lim_{n\to+\infty}x_n-c^*t_n=+\infty$, it follows that the sequences $(U_n)_{n\in\N}$ and $(V_n)_{n\in\N}$ are bounded in $L^{\infty}_{loc}(\R^2)$. From standard parabolic estimates and the fact that $u_{\infty}(t_n,[x_n])\to0$ as $n\to+\infty$, the functions $U_n$ converge locally uniformly in $\R^2$, up to extraction of a subsequence, to a nonnegative classical solution $U_{\infty}$ of
$$(U_{\infty})_t=(U_{\infty})_{xx}+g(x)U_{\infty}\hbox{ in }\R^2$$
(remember that $f'(0)=1$). Furthermore, $(U_n)_x\to(U_{\infty})_x$ locally in $\R^2$ as $n\to+\infty$ and $U_{\infty}(0,0)=1$, whence $U_{\infty}>0$ in $\R^2$ from the maximum principle. In particular, the functions
$$\frac{(u_{\infty})_x(t+t_n,x+[x_n])}{u_{\infty}(t+t_n,x+[x_n])}=\frac{(U_n)_x(t,x)}{U_n(t,x)}$$
are locally bounded. As far as the functions $V_n$ are concerned, they obey
$$\baa{rcl}
(V_n)_t(t,x) & = & \displaystyle(V_n)_{xx}(t,x)+2\frac{(U_n)_x(t,x)}{U_n(t,x)}\,(V_n)_x(t,x)\vspace{3pt}\\
& & \displaystyle+g(x)\,\Big(\frac{f(u_{\infty}(t+t_n,x+[x_n])V_n(t,x))}{u_{\infty}(t+t_n,x+[x_n])}-\frac{f(u_{\infty}(t+t_n,x+[x_n]))}{u_{\infty}(t+t_n,x+[x_n])}V_n(t,x)\Big)\eaa$$
in $\R^2$. Since $(U_n)_x/U_n\to(U_{\infty})_x/U_{\infty}$ and $u_{\infty}(t+t_n,x+[x_n])\to0$ locally uniformly in $\R^2$ as $n\to+\infty$, and since the functions $V_n$ are locally bounded, it follows from standard parabolic estimates that, up to extraction of a subsequence, the functions $V_n$ converge locally uniformly in $\R^2$ to a classical solution $V_{\infty}$ of
\be\label{Vinfty}
(V_{\infty})_t=(V_{\infty})_{xx}+2\frac{(U_{\infty})_x}{U_{\infty}}(V_{\infty})_x\hbox{ in }\R^2.
\ee
Owing to the definitions of $V_n$ and $\xi_*$, one has $V_n\ge 1$ whence $V_{\infty}\ge 1$ in $\R^2$. On the other hand,
$$V_n(\xi_n-\xi_*,x_n-[x_n])=\frac{v^{\xi_n}(t_n,x_n)}{u_{\infty}(t_n,x_n)}\times\frac{U_n(0,x_n-[x_n])}{U_n(\xi_n-\xi_*,x_n-[x_n])}\le\frac{U_n(0,x_n-[x_n])}{U_n(\xi_n-\xi_*,x_n-[x_n])}$$
from~(\ref{xinbis}). By passing to the limit as $n\to+\infty$, one infers that $V_{\infty}(0,x_{\infty})\le 1$. Finally, $V_{\infty}(0,x_{\infty})=1$. Therefore, $V_{\infty}=1$ in $\R^2$ from the strong parabolic maximum principle and the uniqueness of the Cauchy problem associated to~(\ref{Vinfty}).\par
One has then proved that
$$\frac{u_{\infty}(t+t_n+\xi_*,x+[x_n]+1)}{u_{\infty}(t+t_n,x+[x_n])}=\frac{v^{\xi_*}(t+t_n,x+[x_n])}{u_{\infty}(t+t_n,x+[x_n])}\to1\ \hbox{ locally uniformly in }\R^2\hbox{ as }n\to+\infty.$$
It follows by immediate induction that, for each $p\in\N$, there holds
$$\frac{u_{\infty}(t+t_n+p\xi_*,x+[x_n]+p)}{u_{\infty}(t+t_n,x+[x_n])}\to1\ \hbox{ locally uniformly in }\R^2\hbox{ as }n\to+\infty.$$
Fix $p\in\N$. Property~(\ref{exp}) and the limit $\lim_{n\to+\infty}x_n-c^*t_n=+\infty$ imply that, for $n$ large enough,
$$\frac{u_{\infty}(t_n+p\xi_*,[x_n]+p)}{u_{\infty}(t_n,[x_n])}\ge\frac{\kappa\,\big([x_n]+p-c^*t_n-pc^*\xi_*-1\big)\,e^{-\lambda^*([x_n]+p-c^*t_n-pc^*\xi_*)}}{\rho\,([x_n]-c^*t_n)\,e^{-\lambda^*([x_n]-c^*t_n-1)}}.$$
By passing to the limit as $n\to+\infty$, one gets that
$$1\ge\frac{\kappa}{\rho}\,e^{p\lambda^*(c^*\xi_*-1)-\lambda^*}.$$
Since this inequality holds for all $p\in\N$ and since one had assumed that $\xi_*>1/c^*$, one is led to a contradiction. One concludes that $\xi_*\le 1/c^*$, whence $v^{1/c^*}\ge u_{\infty}$ in $\R^2$.\par
By sliding $u_{\infty}(t,x+1)$ in the other $t$-direction, one can prove similarly that $v^{\xi}\le u_{\infty}$ in $\R^2$ for all~$\xi\le\xi_-$ for some real number $\xi_-$, and that the largest such $\xi$ cannot be smaller than $1/c^*$. Therefore, $v^{1/c^*}\le u_{\infty}$ in $\R^2$.\par
Finally, $v^{1/c^*}=u_{\infty}$ in $\R^2$, that is $u_{\infty}(t+1/c^*,x+1)=u_{\infty}(t,x)$ for all $(t,x)\in\R^2$. In other words, $u_{\infty}$ is a pulsating front with speed $c^*$, connecting $0$ and $1$. The conclusion~(\ref{uinftybis}) follows from the uniqueness up to time-shifts of the pulsating fronts, for a given speed (see~\cite{hr}). The proof of Lemma~\ref{lemliouville} is thereby complete.~\hfill$\Box$

%%%%%%%%%%%%%%%%%%%%%%%%%%%%%%%%%%%%%%%%%
%%%%%%%%%%%%%%%%%%%%%%%%%%%%%%%%%%%%%%%%%

\end{document}